\documentclass[dvips,preprint]{imsart}

\RequirePackage[OT1]{fontenc} \RequirePackage{amsthm,amsmath}
\RequirePackage[numbers]{natbib}
\RequirePackage[colorlinks,citecolor=blue,urlcolor=blue]{hyperref}
\RequirePackage{hypernat}
\usepackage{amsmath}
\usepackage{amssymb}
\usepackage{amsmath}
\usepackage{amsfonts}
\usepackage[american]{babel}
\usepackage{graphicx}
\usepackage{flafter}
\usepackage[section]{placeins}
\arxiv{math.PR/0000000}

\startlocaldefs

\def\P{{\mathbb P}}
\def\E{{\mathbb E }}

\def\Bbb E{\mathbb{E}}
\def\Bbb R{\mathbb{R}}
\newtheorem{lemma}{Lemma}
\newtheorem{theorem}{Theorem}

\newtheorem{remark}{Remark}
\makeatletter \@addtoreset{equation}{section}

\makeatother

\font\tencmmib=cmmib10 \skewchar\tencmmib '60
\newfam\cmmibfam
\textfont\cmmibfam=\tencmmib

\font\tenmsb=msbm10 
\def\Bbb#1{\hbox{\tenmsb#1}}

\def\lessim{\ \lower4pt\hbox{$
\buildrel{\displaystyle <}\over\sim$}\ }
\def\gessim{\ \lower4pt\hbox{$\buildrel{\displaystyle >}
\over\sim$}\ }

\def\go0{\to 0}

\def\leftitem#1{\item{\hbox to\parindent{\enspace#1\hfill}}}

\def\sg{\sigma}

\def\sg2{\sigma^2}

\def\__{_{\infty}}

\numberwithin{equation}{section} \theoremstyle{plain}

\newcommand{\1}{{\rm 1}\kern-0.24em{\rm I}}
\def\E{\mathbb E}
\def\R{\mathbb R}
\newtheorem{assumption}{Assumption}

\newtheorem{prop}{Proposition}
\newtheorem{lem}{Lemma}
\newtheorem{defi}{Definition}
\endlocaldefs

\begin{document}

\begin{frontmatter}
\title{Sparse Principal Component Analysis with missing observations} \runtitle{Sparse PCA with missing observations}

\begin{aug}
\author{\fnms{Karim} \snm{Lounici}\thanksref{m1}\ead[label=e2]{klounici@math.gatech.edu}}
\thankstext{m1}{Supported in part by NSF Grant DMS-11-06644 and Simons foundation Grant 209842}
\runauthor{K. Lounici}

\affiliation{Georgia Institute of Technology\thanksmark{m1}}

\address{School of Mathematics\\
Georgia Institute of Technology\\
Atlanta, GA 30332-0160\\
\printead{e2}\\
}
\end{aug}

\begin{abstract}
In this paper, we study the problem of sparse Principal Component
Analysis (PCA) in the high-dimensional setting with missing
observations. Our goal is to estimate the first principal component
when we only have access to partial observations. Existing
estimation techniques are usually derived for fully observed data
sets and require a prior knowledge of the sparsity of the first
principal component in order to achieve good statistical guarantees.
Our contributions is threefold. First, we establish the first
information-theoretic lower bound for the sparse PCA problem with
missing observations. Second, we propose a simple procedure that
does not require any prior knowledge on the sparsity of the unknown
first principal component or any imputation of the missing
observations, adapts to the unknown sparsity of the first principal
component and achieves the optimal rate of estimation up to a
logarithmic factor. Third, if the covariance matrix of interest
admits a sparse first principal component and is in addition
approximately low-rank, then we can derive a completely data-driven
procedure computationally tractable in high-dimension, adaptive to
the unknown sparsity of the first principal component and
statistically optimal (up to a logarithmic factor).

\end{abstract}



\end{frontmatter}

\section{Introduction}

Let $X,X_1,\ldots,X_n \in \R^p$ be i.i.d. zero mean vectors with
unknown covariance matrix $\Sigma = \E X \otimes X$ of the form
\begin{align}\label{model}
\Sigma = \sigma_1 \theta_1\theta_1^\top + \sigma_2 \Upsilon,
\end{align}
where $\sigma_1>\sigma_2 \geq 0$, $\theta_1\in \mathcal S^{p}$ (the $l_2$ unit sphere in $\R^p$) and $\Upsilon$ is a $p\times p$ symmetric
positive semi-definite matrix with spectral norm $\|\Upsilon\|_\infty \leq 1$ and such that $\Upsilon \theta_1 = 0$. The eigenvector
$\theta_1$ is called the first principal component of $\Sigma$. Our objective is to estimate the first principal component $\theta_1$ when the
vectors $X_1,\ldots,X_n$ are partially observed. More precisely, we consider
the following framework. Denote by $X^{(j)}_i$ the $j$\emph{-th}
component of the vector $X_i$. We assume that each component
$X^{(j)}_i$ is observed independently of the others with probability
$\delta\in (0,1]$. Note that $\delta$ can be easily estimated by the
proportion of observed entries. Therefore, we will assume in this
paper that $\delta$ is known. Note also that the case $\delta=1$
corresponds to the standard case of fully observed vectors. Let $(\delta_{i,j})_{1\leq i \leq n,1\leq j \leq p}$ be a
sequence of i.i.d. Bernoulli random variables with parameter
$\delta$ and independent from $X_1,\ldots,X_n$. We observe $n$
i.i.d. random vectors $Y_1,\ldots,Y_n\in \R^p$ whose components
satisfy
\begin{equation}\label{equationY}
Y_i^{(j)} = \delta_{i,j}X_i^{(j)},\quad 1\leq i \leq n,\,1\leq j
\leq p.
\end{equation}
We can think of the $\delta_{i,j}$ as masked variables. If
$\delta_{i,j}=0$, then we cannot observe the $j$\emph{-th} component
of $X_i$  and the default value $0$ is assigned to $Y_i^{(j)}$. Our goal is then to estimate $\theta_1$ given the partial
observations $Y_1,\ldots,Y_n$.

Principal Component Analysis (PCA) is a popular technique to reduce the dimension of a data set that has been used for many years in a variety of different fields including image processing, 
engineering, genetics, meteorology, chemistry and many others. In most of these fields, data are now high-dimensional, that is the number of parameters $p$ is much larger than the sample 
size $n$, and contain missing observations.
This is especially true in genomics with gene expression microarray data where PCA is used to detect the genes responsible for a given biological process. Indeed, despite the
 recent improvments in gene expression techniques, microarray data can contain up to $10\%$ missing observations affecting up to $95
\%$ of the genes. Unfortunately, it is a known fact that PCA is very sensitive even to small perturbations of the data including in particular missing observations. Therefore, several strategies have been
developped to deal with missing values. The simple strategy that consists in eliminating from the PCA study any gene with at least one missing observation is not acceptable in this
 context since up to $95\%$ of the genes can be eliminated from the study. An alternative strategy consists in infering the missing values prior to the PCA using complex imputation schemes \cite{JWWO05,HE09}. 
These schemes usually assume that the genes interactions follow some specified model and involve intensive computational preprocessing to imput the missing observations.
We propose in this paper a different strategy. Instead of building an imputation technique based on assumptions describing the genome structure (about which we usually have
 no prior information), we propose a technique based on the analysis of the perturbations process. In other words, if we understand the process generating the missing observations,
 then we can efficiently correct the data prior to the PCA analysis. This strategy was first introduced in \cite{L12} to estimate the spectrum of low-rank covariance matrices. One of our goal is to show that this approach can be successfully applied
 to perform fast and accurate PCA with missing observations.

Standard PCA in the full observation framework ($\delta=1$) consists in extracting the first principal components of $\Sigma$ (that is the eigenvector $\theta_1$ associated to the
 largest eigenvalue) based on the i.i.d. observations $X_1,\cdots,X_n$:
\begin{align}\label{PCA}
\hat \theta =\mathrm{argmax}_{\theta^\top \theta =1} \theta^\top \Sigma_n \theta,
\end{align}
where $\Sigma_n = \frac{1}{n}\sum_{i=1}^n X_iX_i^\top$. The standard PCA presents two majors drawbacks. First, it is not consistent in high-dimension \cite{JL09,N08,P07}. Second, the solution $\hat\theta$ is usually a dense vector whereas sparse solutions are prefered in most applications in order to obtain simple interpretable structures. For instance, in microarray data,
we typically observe that only a few among the thousands of screened genes are involved in a given biological process. In order to improve interpretability, several approaches have been proposed to perform sparse PCA, that is to enforce sparsity of the PCA outcome. See for instance \cite{PJ07,M11,SH08} for SVD based iterative thresholding approaches.
\cite{ZHT06} reformulated the sparse PCA problem as a sparse regression problem and then used the LASSO estimator. See also \cite{M06} for greedy methods.
We consider now the approach by \cite{Jolliffe} which consists in computing a solution of (\ref{PCA}) under the additional $l_1$-norm constraint $|\theta|_1\leq \bar s$ for some fixed integer $\bar s\geq 1$ in order to enforce sparsity of the solution.
The same approach with the $l_1$-norm constraint replaced by the $l_0$-norm gives the following procedure
\begin{equation}\label{oracle-proc}
 \hat\theta_{\mathrm{o}} = \mathrm{argmax}_{\theta\in\mathcal S^p\,:\, |\theta|_0 \leq \bar s} \left( \theta^\top \Sigma_n \theta  \right),
\end{equation}
where $|\theta|_0$ denotes the number of nonzero components of $\theta$. In a recent paper, \cite{VL12} established
the following oracle inequality
\begin{equation*}
\left(\E\|\hat\theta_{\mathrm{o}}\hat\theta_{\mathrm{o}}^\top -  \theta_1\theta_1^\top\|_2\right)^2 \leq C \left(\frac{\sigma_1}{\sigma_1 -\sigma_2}\right)^2 \bar s \frac{\log(p/\bar s)}{n},
\end{equation*}
for some absolute constant $C>0$. Note that this procedure requires the knowledge of an upper bound $\bar s \geq |\theta_1|_0$.
In practice, we generally do not have access to any prior information on the sparsity of $\theta_1$.
Consequently, if the parameter $\bar s$ we use in the procedure is too small, then the above upper bound does not hold, and if $\bar s$ is too large,
then the above upper bound (even though valid) is sub-optimal. In other words, the procedure (\ref{oracle-proc}) with $\bar s = |\theta_1|_0$
 can be seen as an oracle and our goal is to propose a procedure that performs as well as this oracle without any prior information on $|\theta_1|_0$.

In order to circumvent the fact that $|\theta_1|_0$ is unknown, we consider the following procedure proposed by \cite{ABE08}
\begin{equation}\label{PCA-BIC}
 \hat\theta_1 = \mathrm{argmax}_{\theta\in\mathcal S^p} \left( \theta^\top \Sigma_n \theta -\lambda |\theta|_0 \right),
\end{equation}
where $\lambda>0$ is a regularization parameter to be tuned properly. \cite{ABE08,nest} studied the computational aspect. In particular, \cite{nest} proposed a computationally tractable procedure to solve the above constrained maximization problem even in high-dimension.
However none of these references investigated the statistical performances of this procedure or the question of the optimal tuning of $\lambda$. We propose to carry out this analysis in this paper.
We establish the optimality of this procedure in the minimax sense and explicit the optimal theorethical choice for the regularization parameter $\lambda$.

When the data contains incomplete observations ($\delta<1$), we do not have access to the empirical covariance matrix
$\Sigma_n $. Given the observations $Y_1,\ldots,Y_n$, we can build the following empirical covariance matrix
$$
\Sigma_n^{(\delta)} = \frac{1}{n}\sum_{i=1}^n Y_i Y_i^\top.
$$
As noted in \cite{L12}, $\Sigma_n^{(\delta)}$ is not an unbiased estimator of $\Sigma$, Consequently, we need to consider the following correction in order to get sharp estimation results:
\begin{equation}\label{Sigmaemprecons}
\tilde\Sigma_n = (\delta^{-1} -
\delta^{-2})\mathrm{diag}\left(\Sigma_n^{(\delta)}\right) +
\delta^{-2}\Sigma_n^{(\delta)}.
\end{equation}
Indeed, we can check by elementary algebra that $\tilde \Sigma_n$ is an unbiased estimator of $\Sigma$ in the missing observation framework $\delta\in(0,1]$.
Therefore, we consider the following estimator in the missing observation framework
\begin{equation}\label{esttheta1}
\hat{\theta}_1 = \mathrm{argmax}_{\theta\in\mathcal{ S}^p\,:\, |\theta|_0\leq \bar s} \left(  \theta^\top \tilde \Sigma_n \theta - \lambda |\theta|_0 \right),
\end{equation}
where $\lambda>0$ is a regularization parameter to be tuned properly
and $\bar s$ is a mild constraint on $|\theta_1|_0$. More precisely,
$\bar s$ can be chosen as large as $\frac{\delta^2 n}{\log(ep)}$
when no prior information on $|\theta_1|_0$ is available. We will
prove in particular that the procedure (\ref{esttheta1}) adapts to
the unknown sparsity of $\theta_1$ provided that $|\theta_1|_0 \leq
\bar s$. We also investigate the case where $\Sigma$ is in addition
approximately low-rank. In that case, we can remove the restriction
$|\theta|_0\leq \bar s$ (taking $\bar s =p$) in the procedure
(\ref{esttheta1}) and propose a data-driven choice of the
regularization parameter $\lambda$. We will show that this
data-driven procedure also achieves the optimal rate of estimation
(up to a logarithmic factor) in the missing observation framework
$\delta \in (0,1]$ without any prior knowledge on $|\theta_1|_0$.
Finally, we establish information theoretic lower bounds for the sparse PCA problem in the missing observation framework $\delta\in (0,1]$ with the sharp dependence on $\delta$,
thus expliciting completely the effect of missing observations on
the sparse PCA estimation rate. Note that our results are nonasymptotic in
nature and hold for any setting of $n,p$ including in particular the
high-dimensional setting $p> n$.

The rest of the paper is organized as follows. In Section
\ref{tools}, we recall some tools and definitions that will be
useful for our statistical analysis. Section \ref{secmain} contains
our main theoretical results. 
Finally, Section \ref{proof} contains the proofs of our results.

\section{Tools and definitions}\label{tools}

In this section, we introduce various notations and definitions and we recall some known results that we will use to establish our results.

The $l_q$-norms of a vector
$x=\left(x^{(1)},\cdots,x^{(p)}\right)^{\top}\in \R^p$ is given by
$$
|x|_q = \left(\sum_{j=1}^p |x^{(j)}|^q\right)^{1/q},\; \text{for}\;1
\leq q < \infty, \quad \text{and}\quad |x|_\infty = \max_{1\leq j
\leq p}|x^{(j)}|.
$$
The support of a vector $x=\left(x^{(1)},\cdots,x^{(p)}\right)^{\top}\in \R^p$ is defined as follows
$$
J(x) = \left\lbrace  j\,:\, x^{(j)}\neq 0 \right\rbrace.
$$
We denote the number of nonzero components of $x$ by $|x|_0$. Note that $|x|_0 = |J(x)|$. Set $\mathcal S^p = \left\lbrace x\in \R^p\,:\, |x|_2 =1  \right\rbrace$. For any $J\in [p]$, we define
$\mathcal S^p(J) = \left\{ x\in \mathcal S^p\,:\, J(x) = J \right\} $. For any integer $1\leq s\leq p$, we define $\mathcal S^p_s  = \left\{ x\in \mathcal S^p\,:\, |x|_0 = s  \right\}$.
 Note that $\mathcal S^p_s = \cup_{J\in [p]\,:\, |J| = s} \mathcal S^p(J)$.

For any $p\times p$ symmetric matrix $A$ with eigenvalues $\sigma_1(A), \cdots, \sigma_p(A)$, we define the Schatten $q$-norm of $A$ by
$$
\|A\|_q = \left( \sum_{j=1}^p |\sigma_j(A)|^q\right)^{1/q},\;\forall 1\leq q <\infty,\quad \text{and}\quad \|A\|_{\infty} = \max_{1\leq j \leq p}\left\{|\sigma_j(A)|\right\}.
$$
Define the usual matrix scalar product $\langle A,B\rangle = \mathrm{tr}(A^\top B)$ for any $A,B\in R^{p\times p}$. Note that $\|A\|_2 = \sqrt{\langle A,A\rangle}$ for any $A\in \R^{p\times p}$.
Recall the trace duality property
$$
|\langle A,B\rangle| \leq \|A\|_\infty \|B\|_1,\quad \forall A,B\in \R^{p\times p}.
$$

We recall now some basic facts about $\epsilon$-nets (See for instance Section 5.2.2 in \cite{vershynin}).
\begin{defi}
Let $(A,d)$ be a metric space and let $\epsilon>0$. A subset $\mathcal N_\epsilon$ of $A$ is called an $\epsilon$-net of $A$ if for every point $a\in A$,
 there exists a point $b\in \mathcal N_\epsilon$ so that $d(a,b)\leq \epsilon$.
\end{defi}
We recall now an approximation result of the spectral norm on an $\epsilon$-net.
\begin{lem}\label{eps-net}
Let $A$ be a $k\times k$ symmetric matrix for some $k\geq 1$. For any $\epsilon \in (0,1/2)$, there exists an $\epsilon$-net $\mathcal N_\epsilon \subset \mathcal S^{k}$ (
the unit sphere in $\R^k$) such that
\begin{align*}
|\mathcal N_\epsilon| \leq \left(1+ \frac{2}{\epsilon}\right)^{k},
\end{align*}
and
\begin{align*}
\sup_{\theta\in\mathcal S^k}|\langle Ax,x\rangle | \leq \frac{1}{1-2\epsilon}\sup_{\theta\in \mathcal N_\epsilon}|\langle Ax,x\rangle |.
\end{align*}
\end{lem}
See for instance Lemma 5.2 and Lemma 5.3 in \cite{vershynin} for a proof.

We recall now the definition and some basic properties of
sub-exponential random vectors.
\begin{defi}
The $\psi_\alpha$-norms of a real-valued random variable $V$ are
defined by
$$
\|V\|_{\psi_\alpha} = \inf\left\lbrace u>0: \mathbb E
\exp\left(|V|^\alpha/u^\alpha
 \right)\leq 2 \right\rbrace,\quad \alpha \geq 1.
$$
We say that a random variable $V$ with values in $\R$ is
sub-exponential if $\|V\|_{\psi_\alpha}<\infty$ for some $\alpha
\geq 1$. If $\alpha = 2$, we say that $V$ is sub-gaussian.
\end{defi}
We recall some well-known properties of sub-exponential random
variables:
\begin{enumerate}
\item For any real-valued random variable $V$ such that
$\|V\|_{\alpha}<\infty$ for some $\alpha \geq 1$, we have
\begin{equation}\label{subexp-basic1}
\mathbb E |V|^m \leq 2\frac{m}{\alpha}\Gamma\left(\frac{m}{\alpha}\right)\|V\|_{\psi_{\alpha}}^m,\quad \forall m\geq 1
\end{equation}
where $\Gamma(\cdot)$ is the Gamma function.
\item If a real-valued random variable $V$ is sub-gaussian, then
$V^2$ is sub-exponential. Indeed, we have
\begin{equation}\label{subexp-basic2}
\|V^2\|_{\psi_1}\leq 2 \|V\|_{\psi_2}^2.
\end{equation}
\end{enumerate}

\begin{defi}
A random vector $X\in \R^p$ is sub-exponential if $\langle X,x\rangle$ are sub-exponential random variables for all
$x\in \R^p$. The $\psi_\alpha$-norms of a random vector $X$ are
defined by
$$
\|X\|_{\psi_\alpha} = \sup_{x\in\mathcal S^{p}}\|\langle
X,x\rangle\|_{\psi_\alpha},\quad \alpha\geq 1.
$$
\end{defi}

We recall a version of Bernstein's inequality for unbounded real-valued
random variables.
\begin{prop}\label{Bernstein}
Let $Y_1,\ldots,Y_n$ be independent real-valued random variables with zero mean. Let there exist constants $\sigma$, $\sigma'$ and $K$ such that for any $m\geq 2$
\begin{equation}\label{bernstein-moment}
\frac{1}{n}\sum_{i=1}^n |\E\left[ Y_i^m \right]| \leq \frac{m!}{2}K^{m-2}\sigma^2,\quad \text{and}\quad  \frac{1}{n}\sum_{i=1}^n \E \left[|Y_i|^m \right]\leq \frac{m!}{2}K^{m-2}(\sigma')^2.
\end{equation}
Then for every $t\geq 0$, we have  with probability at least $1-2e^{-t}$
$$
\left| \frac{1}{n}\sum_{i=1}^n  Y_i  \right| \leq \sigma \sqrt{\frac{2t}{n}} + K\frac{t}{n}.
$$
\end{prop}
\begin{remark}
In the usual formulation of Bernstein's inequality, only the second moment condition in (\ref{bernstein-moment}) is imposed and the conclusion holds valid with $\sigma$ replaced by $\sigma'$.
The refinement we propose here is necessary in the missing observation framework in order to get the sharp dependence of our bounds on $\delta$. An investigation of the proof of Bernstein's inequality shows that this refinement follows immediately from Chernoff's bound used to prove the standard Bernstein's inequality (See for instance Proposition 2.9 in \cite{massart}). Indeed, using Chernoff's approach, we need a control on the following expectation
\begin{align*}
\E\left[\Phi(t Y_i)\right] &=  \frac{t^2\E\left[ Y_i^2\right]}{2} + \E\left[\sum_{k=3}^{\infty} \frac{t^k Y_i^k}{k!}\right]= \frac{t^2\E\left[ Y_i^2\right]}{2} + \sum_{k=3}^{\infty} \frac{t^k \E\left[Y_i^k\right]}{k!},\quad \forall t\in \left(0,\frac{1}{K}\right),
\end{align*}
where we have used the second moment condition in (\ref{bernstein-moment}) and Fubini's theorem to justify the inversion of the sum and the expectation. The rest of the proof is left unchanged.
\end{remark}

\section{Main results for sparse PCA with missing observations}\label{secmain}
In this section, we state our main statistical results concerning
the procedure (\ref{esttheta1}). We will establish these results
under the following condition on the distribution of $X$.

\begin{assumption}[Sub-gaussian observations]\label{assumption1}
The random vector $X\in \R^p$ is sub-gaussian, that is $\|
X\|_{\psi_2}<\infty$. In addition, there exist a numerical constant
$c_1>0$ such that
\begin{equation}\label{subexp1}
\E(\langle X, u \rangle)^2\geq c_1 \|\langle X, u
\rangle\|_{\psi_2}^2,\, \forall u\in \R^p.
\end{equation}
\end{assumption}

\subsection{Oracle inequalities for sparse PCA}

We first establish a preliminary result on the stochastic deviation of the following empirical process
\begin{align*}
\mathbf Z_n(s) &= \max_{\theta \in \mathcal S^p_s}\left\{\left| \theta^\top (\tilde \Sigma_n - \Sigma)\theta\right|\right\},\quad \forall 1\leq s \leq p.
\end{align*}
To this end, we introduce the following quantity
\begin{align*}
\zeta_n(s,p,t,\delta) &:= \max\left\lbrace \sqrt{\frac{t+s\log(ep/s)}{\delta^2 n}}, \frac{t+s\log(ep/s)}{\delta^2 n} \right\rbrace.
\end{align*}

\begin{prop}\label{Prop-Zn}
Let Assumption \ref{assumption1} be satisfied. Let $Y_1,\cdots,Y_n$ be defined in (\ref{equationY}) with $\delta\in (0,1]$. Then, we have
\begin{align}
\P\left( \bigcap_{s=1}^p \left\lbrace \mathbf{Z}_n(s) \leq c \frac{\sigma_{\max}(s)}{c_1\wedge 1}\zeta_{n}(s,p,t,\delta)  \right\rbrace \right) \geq 1- e^{-t}
\end{align}
where $c>0$ is an absolute constant and $\sigma_{\max}(s) = \max_{\theta\in\mathcal S^p_s}\left(\theta^\top \Sigma \theta\right)$.
\end{prop}

We can now state our main result
\begin{theorem}\label{main1}
Let Assumption \ref{assumption1} be satisfied. Let $Y_1,\cdots,Y_n$ be defined in (\ref{equationY}) with $\delta\in (0,1]$. Consider the estimator (\ref{esttheta1}) with parameters $\bar s$ satisfying $n\geq \delta^{-2}(\bar s+1)\log(ep/\bar s)$ and
\begin{align}\label{lambda-th}
\lambda = C \frac{\sigma_1^2}{\sigma_1-\sigma_2} \frac{\log (ep)}{\delta^2 n},
\end{align}
where $C>0$ is a large enough numerical constant.\\
If $|\theta_1|_0 \leq \bar s$, then we have, with probability at least $1-\frac{1}{p}$, that
\begin{align*}
\|\hat\theta_1\hat\theta_1^\top - \theta_1\theta_1^\top \|_2^2 \leq C'|\theta_1|_0 \tilde{\sigma}^2\frac{\log(ep)}{\delta^2 n}.
\end{align*}
where $\tilde \sigma =  \frac{\sigma_1}{\sigma_1-\sigma_2}$ and $C'>0$ is a numerical constant that can depend only on $c_1$.
\end{theorem}
\begin{enumerate}
\item We observe that the estimation bound increases as the difference $\sigma_1-\sigma_2$ decreases. The problem of estimation of the first principal component is statistically more difficult when the largest and second largest eigenvalues are close.
We also observe that the optimal choice of the regularization
parameter (\ref{lambda-th}) depends on the eigenvalues
$\sigma_1,\sigma_2$ of $\Sigma$. Unfortunately, these quantities are
typically unknown in practice. In order to circumvent this
difficulty, we propose in Section \ref{datadriven} a data-driven
choice of $\lambda$ with optimal statistical performances (up to a
logarithmic factor) provided that $\Sigma$ is approximately
low-rank.

\item Let now consider the full observation framework ($\delta=1$). In that case, if $|\theta_1|_0\leq \bar s$, we obtain the following upper bound with probability at least $1-\frac{1}{p}$
\begin{align*}
\|\hat\theta_1\hat\theta_1^\top - \theta_1\theta_1^\top \|_2^2 \leq C'|\theta_1|_0 \tilde{\sigma}^2\frac{\log(ep)}{n}.
\end{align*}
We can compare this result with that obtained  for the procedure
(\ref{oracle-proc}) by \cite{VL12}
\begin{align*}
\left(\E\|\hat\theta_o\hat\theta_o^\top - \theta_1\theta_1^\top \|_2\right)^2 \leq C'\bar s \tilde{\sigma}^2\frac{\log(ep/\bar s)}{n}.
\end{align*}
We see that in order to achieve the rate $|\theta_1|_0\tilde
\sigma^2 \log(ep/|\theta_1|_0)$ with the procedure
(\ref{oracle-proc}), we need to know the sparsity of $\theta_1$ in
advance, whereas our procedure adapts to the unknown sparsity of
$\theta_1$ and achieves the minimax optimal rate up to a logarithmic
factor provided that $|\theta_1|_0\leq \bar s$ (see Section
\ref{lower} for the lower bounds). This logarithmic factor is the
price we pay for adaptation to the sparsity of $\theta_1$. Note also that we can formulate a version of (\ref{oracle-proc}) when observations are missing ($\delta<1$)
 by replacing $\Sigma_n$ with $\tilde\Sigma_n$. In that case, our techniques of proof will give with probability at least $1-\frac{1}{p}$
\begin{align*}
\|\hat\theta_{\mathrm{o}}\hat\theta_{\mathrm{o}}^\top - \theta_1\theta_1^\top \|_2^2 \leq C'\bar s \tilde{\sigma}^2\frac{\log(ep/\bar s)}{\delta^2 n}.
\end{align*}

\item We discuss now the choice of $\bar s$. In practice, when no prior information
on the sparsity of $\theta_1$ is available, we propose to choose
$\bar s = \frac{\delta^2 n}{\log(ep)}-1$. Then, the procedure
(\ref{esttheta1}) adapts to the unknown sparsity of $\theta_1$
provided that the condition $|\theta_1|_0 \leq \frac{\delta^2
n}{\log(ep)}-1$ is satisfied, which is actually a natural condition
on $\theta_1$ in order to obtain a non trivial estimation result.
Indeed, if $|\theta_1|_0 \geq \delta^2 n/\log( ep)$, then the upper
bound in Theorem \ref{main1} for the estimator $\hat\theta_1$
becomes larger than $\tilde \sigma^2>1$ whereas the bound for the
null estimator is $\|0-\theta_1\theta_1^\top \|_2^2 =1$.

\item In the case where observations are missing ($\delta<1$), Theorem \ref{main1} guarantees that recovery
of the first principal component is still possible using the
procedure (\ref{esttheta1}). We observe the additional factor
$\delta^{-2}$. Consequently, the estimation accuracy of the
procedure (\ref{esttheta1}) will decrease as the proportion of
observed entries $\delta$ decreases. We show in Section \ref{lower}
below that the dependence of our bounds on $\delta^{-2}$ is sharp.
In other words, there exists no statistical procedure that achieves
an upper bound without the factor $\delta^{-2}$. Thus, we can
conclude that the factor $\delta^{-2}$ is the statistical price to
pay to deal with missing observations in the principal component
estimation problem. If we consider for instance microarray datasets
where typically about $10\%$ of the observations are missing (that
is $\delta=0.9$), then the optimal bound achieved for the first
principal component estimation increases by a factor $1.24$ as
compared to the full observation framework ($\delta=1$).
\end{enumerate}

\subsection{Study of approximately low-rank covariance matrices}

We now assume that $\Sigma$ defined in $(\ref{model})$ is also
approximately low-rank and study the different implications of this
additional condition. We recall that the effective rank of $\Sigma$
is defined by $\mathbf{r}(\Sigma) =
\mathrm{tr}(\Sigma)/\|\Sigma\|_\infty$ where $\mathrm{tr}(\Sigma)$
is the trace of $\Sigma$. We say that $\Sigma$ is approximately
low-rank when $\mathbf r(\Sigma)\ll p$. Note also that the effective
rank of a covariance matrix can be estimated efficiently by
$\mathbf{r}(\hat\Sigma)$ where $\hat\Sigma$ is an acceptable estimator of
$\Sigma$. See \cite{L12} for more details. Thus, the approximately
low-rank assumption can easily be checked in practice.

First, we can propose a different control of the stochastic
quantities $\mathbf Z_n(s)$. Note indeed that $\mathbf Z_n(s) \leq \|\tilde
\Sigma_n - \Sigma\|_\infty$ for any $1\leq s \leq p$. We apply now
Proposition 3 in \cite{L12} and get the following control on
$\mathbf Z_n(s)$. Under the assumptions of Proposition
\ref{Prop-Zn}, we have with probability at least $1-e^{-t}$ that
\begin{equation}\label{spectrum-bound}
\|\tilde\Sigma_n - \Sigma\|_\infty \leq
C\frac{\sigma_1}{c_1} \max\left\{
\sqrt{\frac{\mathbf{r}(\Sigma)\left(t+ \log(2p)\right)}{\delta^2
n}},\frac{\mathbf{r}(\Sigma)\left(t+ \log(2p)\right)}{\delta^2
n}\left(c_1\delta + t+ \log n\right)\right\},
\end{equation}
where $C>0$ is an absolute constant. We concentrate now on the high-dimensional setting $p>n$. Assume that
\begin{align}\label{nlowrank}
n\geq c\frac{\mathbf{r}(\Sigma) \log^2(ep)}{\delta^2},
\end{align}
for some sufficiently large numerical constant $c>0$.
Taking $t=\log(ep)$, we get from the two above displays, with probability at least $1-\frac{1}{ep}$ that
\begin{equation*}
\|\tilde\Sigma_n - \Sigma\|_\infty \leq
c'\sigma_1 \sqrt{\frac{\mathbf{r}(\Sigma)\log(ep)}{\delta^2
n}},
\end{equation*}
where $c'>0$ can depends only on $c_1$. Combining the previous display with Proposition \ref{Prop-Zn} and a union bound argument, we immediately obtain the following control on $\mathbf{Z}_n(s)$.
\begin{prop}\label{prop-Zn-2}
Let the conditions of Proposition \ref{Prop-Zn} be satisfied. In addition, let (\ref{nlowrank}) be satisfied. Then we have
\begin{equation*}
\P\left( \bigcap_{s=1}^p \left\lbrace \mathbf{Z}_n(s) \leq
c\sigma_1 \sqrt{ \min\{s ,\mathbf{r}(\Sigma)\}} \sqrt{\frac{\log(ep)}{\delta^2
n}}\right\rbrace \right)\geq 1-\frac{1}{p},
\end{equation*}
where $c>0$ is numerical constant that can depend only on $c_1$.
\end{prop}
The motivation behind this new bound is the following. When
(\ref{nlowrank}) is satisfied, we can remove the restriction
$|\theta|_0 \leq \bar s$ in the procedure (\ref{esttheta1}). We
consider now
\begin{equation}\label{esttheta2}
\tilde{\theta}_1 = \mathrm{argmax}_{\theta\in\mathcal{ S}^p} \left(  \theta^\top \tilde \Sigma_n \theta - \lambda |\theta|_0 \right).
\end{equation}
Then, a solution of this problem can be computed efficiently even in
the high-dimensional setting using a generalized power method (see
\cite{nest} for more details on the computational aspect), whereas it
is not clear whether the same holds true for the procedure
(\ref{esttheta1}) with the constraint $|\theta_1|_0\leq \bar s$.

We now consider the statistical performance of the procedure
(\ref{esttheta2}). Following the proof of Theorem \ref{main1}, we
establish this result for $\tilde\theta_1$.
\begin{theorem}\label{main2}
Let Assumption \ref{assumption1} be satisfied. Let $Y_1,\cdots,Y_n$ be defined in (\ref{equationY}) with $\delta\in (0,1]$. In addition, let (\ref{nlowrank}) be satisfied. Take
\begin{align}\label{lambda-th}
\lambda = C \frac{\sigma_1^2}{\sigma_1-\sigma_2} \frac{\log (ep)}{\delta^2 n},
\end{align}
where $C>0$ is a large enough numerical constant. Then we have, with probability at least $1-\frac{1}{p}$, that
\begin{align*}
\|\tilde\theta_1\tilde\theta_1^\top - \theta_1\theta_1^\top \|_2^2 \leq C'|\theta_1|_0 \tilde{\sigma}^2\frac{\log(ep)}{\delta^2 n}.
\end{align*}
where $\tilde \sigma =  \frac{\sigma_1}{\sigma_1-\sigma_2}$ and $C'>0$ is a numerical constant that can depend only on $c_1$.
\end{theorem}
Note that this result holds without any condition on the sparsity of
$\theta_1$. Of course, as we already commented for Theorem
\ref{main1}, the result is of statistical interest only when
$\theta_1$ is sparse: $|\theta_1|_0 \leq  \frac{\delta^2 n}{\tilde
\sigma^2 \log (ep)}$. The interest of this result is to guarantee
that the computationally tractable estimator (\ref{esttheta2}) is
also statistically optimal.

\subsection{Data-driven choice of $\lambda$}\label{datadriven}
As we see in Theorem \ref{main2}, the optimal choice of the
regularization parameter depends on the largest and second largest
eigenvalues of $\Sigma$.
 These quantities are typically unknown in practice. To circumvent this difficulty,
  we propose the following data-driven choice for the regularization parameter $\lambda$
\begin{align}\label{lambdadriven}
\lambda_{D} = C\frac{\hat\sigma_1^2}{\hat\sigma_1 - \hat\sigma_2}\frac{\log(ep)}{\delta^2n},
\end{align}
where $C>0$ is a numerical constant and $\hat\sigma_1$ and $\hat \sigma_2$ are the two largest eigenvalues of $\tilde \Sigma_n$.
If (\ref{nlowrank}) is satisfied, then as a consequence of Proposition 3 in \cite{L12}, $\hat\sigma_1$ and $\hat \sigma_2$ are good estimators
 of $\sigma_1$ and $\sigma_2$ even in the missing observation framework. In order to guarantee that $\lambda_D$ is a suitable choice,
we will need a more restrictive condition on the number of
measurements $n$ than (\ref{nlowrank}). This new condition involves
in addition the "variance" $\tilde \sigma^2 =\frac{\sigma_1^2}{(\sigma_1 -\sigma_2)^2}$:
\begin{align}\label{nlowrank2}
n \geq c \frac{\tilde
\sigma^2}{\delta^2}\mathbf{r}(\Sigma)\log^2(ep),
\end{align}
where $c>0$ is a sufficiently large numerical constant. As
compared to (\ref{nlowrank}), we observe the additional factor
$\tilde \sigma^2$ in the above condition. We already noted that
matrices $\Sigma$ for which the difference $\sigma_1-\sigma_2$ is
small are statistically more difficult to estimate. We observe that
the number of measurements needed to construct a suitable
data-driven estimator also increases as the difference
$\sigma_1-\sigma_2$ decreases to $0$.

We have the following result.
\begin{lemma}\label{lem-datadriven}
Let the conditions of Proposition \ref{Prop-Zn} be satisfied, Assume in addition that (\ref{nlowrank2}) is satisfied. Let $\lambda_D$ be defined in (\ref{lambdadriven}). Then, we have with probability at least $1-\frac{1}{p}$ that
$$
\mathbf Z_n(s) \leq \lambda_D\, s,\quad \forall 1\leq s \leq p,
$$
and
$$
\lambda_{D} \leq  C'\frac{\sigma_1^2}{\sigma_1 - \sigma_2}\frac{\log(ep)}{\delta^2n},
$$
for some numerical constant $C'>0$.
\end{lemma}

Consequently, the conclusion of Theorem \ref{main2} holds true for the estimator
(\ref{esttheta2}) with $\lambda = \lambda_D$ provided that (\ref{nlowrank2}) is satisfied.

\subsection{Information theoretic lower bounds}\label{lower}
We derive now minimax lower bounds for the estimation of the first principal component $\theta_1$ in the missing observation framework.

Let $s_1\geq 1$. We denote by $\mathcal{C} = \mathcal
C_{s_1}(\sigma_1,\sigma_2)$ the class of covariance matrices
$\Sigma$ satisfying (\ref{model}) with $\sigma_1>\sigma_2$,
$\theta_1\in \mathcal S^{p}$ with $|\theta|_0\leq s_1$ and
$\Upsilon$ is a $p\times p$ symmetric positive semi-definite matrix
with spectral norm $\|\Upsilon\|_\infty \leq 1$ and such that
$\Upsilon \theta_1 = 0$. We prove now that the dependence of our
estimation bounds on $\sigma_1-\sigma_2,\delta, s_1,n,p$ in Theorems
\ref{main1} and \ref{main2} is sharp in the minimax sense. Set $\bar
\sigma^2 = \frac{\sigma_1\sigma_2}{(\sigma_1 - \sigma_2)^2}$.

\begin{theorem}\label{thlowerbound}
Fix $\delta \in (0,1]$ and $s_1 \geq 3$. Let the integers $n,p\geq
3$ satisfy
\begin{align}\label{nthreshold}
2\bar \sigma^2 s_1\log(ep/s_1)\leq \delta^{2} n.
\end{align}
Let $X_1,\ldots,X_n$ be i.i.d. random vectors in $\R^p$ with
covariance matrix $\Sigma \in \mathcal C$. We observe $n$ i.i.d.
random vectors $Y_1,\ldots,Y_n\in\R^p$ such that
$$
Y_{i}^{j} = \delta_{i,j}X_i^{(j)},\; 1\leq i \leq n,\; 1\leq j \leq
p,
$$
where $(\delta_{i,j})_{1\leq i\leq n,\,1\leq j \leq p}$ is an i.i.d.
sequence of Bernoulli $B(\delta)$ random variables independent of
$X_1,\ldots,X_n$.

Then, there exist absolute constants $\beta\in(0,1)$ and $c>0$ such that
\begin{equation}\label{eq:lower1}
\inf_{\hat{\theta}_1}
\sup_{\substack{\Sigma\in\,{\cal C}
}}
\P_\Sigma\bigg(\|\hat\theta_1\hat\theta_1^\top-\theta_1\theta_1^\top\|^2_{2}>
c \bar\sigma^2\frac{s_1}{\delta^2 n} \log\left(\frac{ep}{s_1}\right)
\bigg)\ \geq\ \beta,
\end{equation}
where $\inf_{\hat{\theta}_1}$ denotes the infimum over all possible
estimators $\hat{\theta}_1$ of $\theta_1$ based on $Y_1,\ldots,Y_n$.
\end{theorem}

\begin{remark}
For $s_1 =1$, we can prove a similar lower bound with the factor
$\delta^{-2}$ replaced by $\delta^{-1}$. This is actually the right
dependence on $\delta$ for $1$-sparse vectors. We can indeed derive
an upper bound of the same order for the selector $e_{\hat\jmath} =
\mathrm{argmax}_{1\leq j \leq p}\left( e_j^\top \tilde \Sigma_n e_j
\right)$ where $e_1,\ldots,e_p$ are the canonical vectors of $\R^p$.\\
For $s_1=2$, we can prove a lower bound of the form
(\ref{eq:lower1}) without the logarithmic factor by comparing for
instance the hypothesis $\theta_0 = \frac{1}{\sqrt{2}}(e_1 + e_2)$
and $\theta_1 =\frac{1}{2}e_1 + \frac{\sqrt{3}}{2}e_2$. Getting a
lower bound for $s_1=2$ with the logarithmic factor remains an open
question.
\end{remark}

\section{Proofs}\label{proof}

\subsection{Proof of Proposition \ref{Prop-Zn}}

\begin{proof}
For any $s\geq1$, we have
\begin{align}\label{proc-Z}
\mathbf Z_{n}(s) &\leq  \delta^{-1}\mathbf Z_n^{(1)}(s) +\delta^{-2}\mathbf Z_n^{(2)}(s)
\end{align}
where
\begin{align*}
\mathbf Z_n^{(1)}(s) =\max_{\theta\in \mathcal S^p_s} \left\lbrace\left| \theta^\top \mathrm{diag}\left(\Sigma_n^{(\delta)} -  \Sigma^{(\delta)}\right)\theta\right|\right\rbrace, \quad \mathbf Z_n^{(2)}(s) =\max_{\theta\in \mathcal S^p_s} \left\lbrace\left| \theta^\top \left(A_n^{(\delta)} -A^{(\delta)}\right)\theta\right|\right\rbrace
\end{align*}
with $A_n^{(\delta)} = \Sigma^{(\delta)}_n - \mathrm{diag}(\Sigma^{(\delta)}_n) $ and $ A^{(\delta)} = \Sigma^{(\delta)} - \mathrm{diag}(\Sigma^{(\delta)})$.

Before we proceed with the study of the empirical processes $\mathbf Z_n^{(1)}(s)$ and $\mathbf Z_n^{(2)}(s)$, we need to introduce first some additional notations. Define
$$Y=(\delta_{1}X^{(1)},\ldots,\delta_{p}X^{(p)})^{\top},$$
where $\delta_{1},\ldots,\delta_{p}$ are i.i.d. Bernoulli random variables
with parameter $\delta$ and independent from $X$. Denote by $\mathbb E_{\delta}$ and $\mathbb E_X$ the expectations w.r.t. $(\delta_1,\cdots,\delta_p)$ and $X$ respectively. We also note that for any $m\geq 2$ and any $\theta\in\mathcal S^p$,
we have $\E\left[\left(\theta^\top YY^\top \theta\right)^m\right] = \E_X\E_\delta\left[\left(\theta^\top YY^\top \theta\right)^m\right] = \E_\delta\E_X\left[\left(\theta^\top YY^\top \theta\right)^m\right]$.
This comes from Fubini's theorem and the fact that $X$ is sub-gaussian.

We now proceed with the study of $\mathbf Z_n^{(2)}(s)$. For any $s\geq 1$ and any fixed $\theta\in \mathcal S^p_s$, we have
\begin{align*}
\theta^\top \left( A_n^{(\delta)} - A^{(\delta)}  \right)\theta &= \frac{1}{n}\sum_{i=1}^n \left[ \theta^\top \left( Y_iY_i^\top - \mathrm{diag}(Y_iY_i^\top) \right) \theta  - \delta^2\theta^\top\left(\Sigma - \mathrm{diag}(\Sigma)\right)\theta \right].
\end{align*}
Set $$Z_i =\left[ \left( Y_iY_i^\top -
\mathrm{diag}(Y_iY_i^\top) \right)  -
\delta^2\left(\Sigma - \mathrm{diag}(\Sigma)\right)
\right],$$
and
$$Z =\left[ \left( YY^\top -
\mathrm{diag}(YY^\top) \right)  -
\delta^2\left(\Sigma - \mathrm{diag}(\Sigma)\right)
\right].$$ We note that $Z,Z_1,\cdots,Z_n$ are i.i.d. We now study the moments $\E\left[\left( \theta^\top Z \theta \right)^m\right]$
 and $\E\left[\left| \theta^\top Z \theta \right|^m\right]$ for any $\theta\in \mathcal{S}^p_s$ and $m\geq 2$.

Note first that
$$
|\theta^{\top}(\Sigma-\mathrm{diag}(\Sigma))\theta|
\leq\max\left(\theta^{\top}\Sigma\theta\, ,\,\theta^{\top}\mathrm{diag}(\Sigma)\theta\right)\leq \sigma_{\max}(s),\quad
\forall \theta\in \mathcal{S}^p_s,\,\forall s\geq 1,
$$
where $\sigma_{\max}(s) = \max_{u\in \mathcal S^p_s}\left(u^\top \Sigma u\right)$. (We have indeed that $\max_{u\in \mathcal S^p_s}\left(u^{\top}\mathrm{diag}(\Sigma)u\right) = \max_{1\leq j \leq p}(\Sigma_{j,j})\leq \sigma_{\max}(s)$.

Thus we get, for any $m\geq 2$ and $\theta\in\mathcal S^p_s$, that
\begin{align}\label{condmoment1}
\E\left[\left(\theta^\top Z\theta\right)^m \right]&\leq 2^{m-1}\E \left[ \left(\theta^\top \left( YY^\top - \mathrm{diag}(YY^\top) \right) \theta \right)^m\right] + 2^{m-1}\delta^{2m}\sigma_{\max}^m(s)
\end{align}
and
\begin{align}\label{condmoment2}
\E\left[\left|\theta^\top Z\theta\right|^m\right] &\leq 2^{m-1}\E \left[ \left|\theta^\top \left( YY^\top - \mathrm{diag}(YY^\top) \right) \theta \right|^m\right] + 2^{m-1}\delta^{2m}\sigma_{\max}^m(s)
\end{align}

For any $\theta\in S^p$ and $\delta_1,\cdots,\delta_p$, we set
$\theta_\delta = (\delta_1\theta^{(1)},\cdots,\delta_p\theta^{(p)}
)^\top$ and we note the following simple fact that we will use in
the next display: $|\theta_\delta|_0 \leq |\theta|_0$ with
probability $1$. We also set $W = XX^\top +(\delta-1) \mathrm{diag}(XX^\top) = (w_{j,k})_{1\leq j,k\leq p}$.
 Note that $w_{j,k} = X^{(j)}X^{(k)}$ for any $j\neq k$ and $w_{j,j} = \delta \left(X^{(j)}\right)^2$ for $j=1,\ldots,p$.

For any $m\geq 2$ and $\theta\in\mathcal S^p_s$, we have
\begin{align}\label{moment-1}
\E\left[ \left(\theta^\top (YY^\top -
\mathrm{diag}(YY^\top))\theta\right)^m \right]&= \E\left[ \left(\theta_\delta^\top (XX^\top -
\mathrm{diag}(XX^\top))\theta_\delta\right)^m\right] \notag\\
&\leq 2^{m-1} \E \left[ \left( \theta_\delta^\top W \theta_\delta \right)^{m}  \right]\notag\\ &\hspace{2cm}+ 2^{m-1}\delta^{m} \E \left[ \left( \theta_\delta^\top \mathrm{diag}(XX^\top) \theta_\delta \right)^{m}\right]
\end{align}
Next, we have for any $m\geq 2$ and $\theta\in\mathcal S^p_s$ that
\begin{align*}
\left(\theta_\delta^\top W\theta_\delta\right)^m = \sum_{j_1,k_1=1}^{p}\cdots\sum_{j_m,k_m=1}^{p} \prod_{t=1}^{m}\theta_{\delta}^{(j_t)}\theta_\delta^{(k_t)}w_{j_t,k_t}.
\end{align*}
Taking the expectation, we get for any $m\geq 2$ and $\theta\in\mathcal S^p_s$ that
\begin{align}\label{moment-2}
\E\left(\theta_\delta^\top W\theta_\delta\right)^m &=\E_X\E_\delta \left[\sum_{j_1,k_1=1}^{p}\cdots\sum_{j_m,k_m=1}^{p} \prod_{t=1}^{m}\theta_{\delta}^{(j_t)}\theta_\delta^{(k_t)}w_{j_t,k_t}\right]\notag\\
&=\E_X \left[\delta^{m} \sum_{j_1,k_1=1}^{p}\cdots\sum_{j_m,k_m=1}^{p} \prod_{t=1}^{m}\theta^{(j_t)}\theta^{(k_t)}X^{(j_t)}X^{(k_t)} \right]\notag\\
&= \delta^{m} \E_X\left[ (\theta^\top X)^{2m} \right]\notag\\
&\leq 2\delta^m m!\left(2\|\theta^\top X\|_{\psi_2}^2 \right)^{m}\notag\\
&\leq 2\delta^m m!\left(\frac{2}{c_1}\sigma_{\max}(s)\right)^{m},
\end{align}
where we have used (\ref{subexp-basic1}), (\ref{subexp-basic2}) and Assumption \ref{assumption1} in the last two lines.

Similarly, we have for any $m\geq 2$ and $\theta\in\mathcal S^p_s$ that
\begin{align*}
\left( \theta_\delta^\top \mathrm{diag}(XX^\top) \theta_\delta \right)^{m} &= \left(\sum_{j=1}^{p} \left(\theta_\delta^{(j)} X^{(j)}\right)^2  \right)^m\\
&\leq \left(\sum_{j=1}^{p} \left(\theta^{(j)} X^{(j)}\right)^2  \right)^m\\
&\leq \sum_{j=1}^p \left(\theta^{(j)}\right)^2 \left( X^{(j)} \right)^{2m},
\end{align*}
where we have used the convexity of $x\rightarrow x^m$ and the fact that $\theta\in \mathcal S^p$ in the last line. Taking now the expectation, we get for any $m\geq 2$ and $\theta\in\mathcal S^p_s$
\begin{align}\label{moment-3}
\E\left( \theta_\delta^\top \mathrm{diag}(XX^\top) \theta_\delta \right)^{m}& \leq \sum_{j=1}^p \left(\theta^{(j)}\right)^2 \E\left[\left( X^{(j)} \right)^{2m}\right]\notag\\
&\leq 2m! \left( 2\max_{1\leq j \leq p}\|X^{(j)}\|_{\psi_2}^2 \right)^m\notag\\
&\leq 2m! \left( \frac{2}{c_1}\max_{j}\left(\Sigma_{j,j}\right) \right)^{m}\notag\\
&\leq 2m! \left(\frac{2}{c_1}\sigma_{\max}(s) \right)^m,
\end{align}
where we have used again (\ref{subexp-basic1}) and (\ref{subexp-basic2}) in the second line line and Assumption \ref{assumption1} in the third line.

Combining (\ref{moment-1}), (\ref{moment-2}) and (\ref{moment-3}), we get for any $m\geq 2$ and any $\theta\in \mathcal S^p_s$ that
\begin{align}
\E\left[ \left(\theta^\top (YY^\top -
\mathrm{diag}(YY^\top))\theta\right)^m \right] &\leq 2m!\left(\frac{ 4\delta}{c_1} \sigma_{\max}(s)\right)^m.
\end{align}
We now plug the above bound in (\ref{condmoment1}) to get for any $m\geq 2$ and any $\theta\in \mathcal S^p_s$ that
\begin{align}\label{condmoment3}
\E\left[ \left(\theta^\top Z \theta\right)^m \right] &\leq \frac{m!}{2}\left(\frac{ 8\sqrt{2}\delta}{c_1} \sigma_{\max}(s)\right)^2 \left(\frac{8\sigma_{\max}(s)}{c_1}\right)^{m-2}+\frac{1}{2} \left( 2\delta^2 \sigma_{\max}(s) \right)^m\notag\\
&\leq\frac{m!}{2}\left(\frac{ 16\delta}{c_1\wedge 1} \sigma_{\max}(s)\right)^2 \left(\frac{8\sigma_{\max}(s)}{c_1\wedge 1}\right)^{m-2}.
\end{align}

Similar (and actually faster computations) give for any $m\geq 2$ and any $\theta\in \mathcal S^p_s$ that
\begin{align*}
&\E\left[ \left| \theta^\top (YY^\top - \mathrm{diag}(YY^\top))\theta \right|^m \right]\\
&\hspace{4cm}\leq 2^{m-1}\E_\delta \E_X \left[ \left(\theta_\delta^\top X\right)^{2m}\right] +  2^{m-1} \E_\delta \E_X \left[ \left( \theta_\delta^\top \mathrm{diag}(XX^\top) \theta_\delta \right)^{m}\right]\\
&\hspace{4cm}\leq 2^{m}m! \E_{\delta}\left[ \left( 2\|\theta_\delta^\top X\|_{\psi_2}^2\right)^{m} \right] + 2^m m! \left(\frac{2}{c_1}\sigma_{\max}(s)\right)^{m}\\
&\hspace{4cm}\leq 2^{m}m! \left( 2\|\theta^\top X\|_{\psi_2}^2\right)^{m}+  m! \left(\frac{4}{c_1}\sigma_{\max}(s)\right)^{m}\\
&\hspace{4cm}\leq 2 m! \left( \frac{4}{c_1}\sigma_{\max}(s) \right)^m,
\end{align*}
since for any $\delta = (\delta_1,\ldots,\delta_p)\in \{0,1\}^p$, $\|\theta_\delta^\top X\|_{\psi_2}\leq \|\theta^\top X\|_{\psi_2}$.
Thus, we get for any $m\geq 2$ and any $\theta\in \mathcal S^p_s$ that
\begin{align}\label{condmoment4}
\mathbb E\left[ \left|\theta^\top Z \theta\right|^m \right] &\leq \frac{m!}{2}\left(\frac{16}{c_{1} \wedge 1} \sigma_{\max}(s)\right)^2 \left(\frac{8}{c_{1}\wedge 1}\sigma_{\max}(s)\right)^{m-2}.
\end{align}

We see that the moments conditions in (\ref{bernstein-moment}) are satisfied with $K =\frac{8}{c_1\wedge 1}\sigma_{\max}(s)$, $\sigma' = \frac{ 16}{c_1\wedge 1} \sigma_{\max}(s)$ and $\sigma = \delta \sigma'$. Thus, for any fixed $\theta\in \mathcal S^p_s$, Bernstein's
inequality gives for any $t'>0$ that
\begin{align*}
\P\left(\left|\theta^\top (A_n^{(\delta)} - A^{(\delta)})\theta
\right| > \frac{16\sqrt{2}}{c_1\wedge 1}\delta \sigma_{\max}(s)\sqrt{\frac{t'}{n}} + \frac{8\sigma_{\max}(s)}{c_1\wedge 1} \frac{t'}{n} \right) \leq 2e^{-t'}
\end{align*}
Note now that
\begin{align*}
\mathbf Z_n^{(2)}(s) = \max_{\theta\in \mathcal S^p_s} \left\lbrace\left| \theta^\top (A_n^{(\delta)} - A^{(\delta)})\theta\right|\right\rbrace = \max_{J\in [p]\,:\, |J|=s}\;\max_{\theta\in \mathcal S^p(J)} \left\lbrace\left| \theta^\top (A_n^{(\delta)} - A^{(\delta)})\theta\right|\right\rbrace.
\end{align*}
For any fixed $J\in [p]$ such that $|J|=s$, Lemma \ref{eps-net} guarantees the existence of a $\frac{1}{4}$-net $\mathcal N(J)$ sucht that $|\mathcal N(J)|\leq 9^s$ and
\begin{align*}
\max_{\theta\in \mathcal S^p(J)} \left\lbrace\left| \theta^\top (A_n^{(\delta)} - A^{(\delta)})\theta\right|\right\rbrace \leq 2 \max_{\theta\in \mathcal N(J)} \left\lbrace\left| \theta^\top (A_n^{(\delta)} - A^{(\delta)})\theta\right|\right\rbrace.
\end{align*}
Combining the last three displays with an union bound argument, we get for
$t' = t + s\log( 9) +  s\log\left(\frac{ep}{s}\right)$ and $t>0$ that
\begin{align}\label{bernstein-Z2}
\P\left( \mathbf Z_n^{(2)}(s) > \zeta_n^{(2)}(s,t)  \right)\leq 2e^{-t},
\end{align}
with
\begin{align*}
\zeta_n^{(2)}(s,t)= \frac{8\sigma_{\max}(s)}{c_1\wedge 1}\left(2\sqrt{2}\delta \sqrt{\frac{t + s\log( 9) +  s\log\left(\frac{ep}{s}\right)}{n}} +  \frac{t + s\log( 9) +  s\log\left(\frac{ep}{s}\right)}{n}\right)
\end{align*}

We proceed similarly to treat the quantity $\mathbf Z_n^{(1)}(s)$. We first note that
\begin{align*}
\theta^\top \left( \mathrm{diag}(\Sigma^{(\delta)}_n - \Sigma^{(\delta)})\right)\theta &= \frac{1}{n}\sum_{i=1}^n \sum_{j=1}^p \left(\left[\theta^{(j)} Y_i^{(j)}\right]^2 - \delta \Sigma_{j,j}\left(\theta^{(j)}\right)^2 \right).
\end{align*}
Next, proceeding essentially as in (\ref{moment-3}), we get for any $m\geq 2$ and any $\theta\in\mathcal S^p_s$
\begin{align*}
&\E\left[\left(\sum_{j=1}^p \left(\theta^{(j)}Y^{(j)}\right)^2 - \delta \Sigma_{j,j}\left(\theta^{(j)}\right)^2  \right)^m\right]\\
&\hspace{4cm}\leq \sum_{j=1}^p\left(\theta^{(j)}\right)^{2}\left(  \E\left[ \left(Y^{(j)}\right)^{2m}\right] +\delta^m \Sigma_{j,j}^m \right)\\
&\hspace{4cm}\leq   \sum_{j=1}^p(\theta^{(j)})^{2}\left[ \delta \E_X\left[ \left(X^{(j)}\right)^{2m}\right] +\delta \Sigma_{j,j}^m \right]\\
&\hspace{4cm}\leq \sum_{j=1}^p(\theta^{(j)})^{2}\left[ 2\delta m! \left(2\|X^{(j)}\|_{\psi_2}^2\right)^m +\delta^m \Sigma_{j,j}^m \right]\\
&\hspace{4cm}\leq \sum_{j=1}^p(\theta^{(j)})^{2}\left[ 2\delta m! \left(\frac{2}{c_1}\Sigma_{j,j}\right)^m +\delta^m \Sigma_{j,j}^m \right]\\
&\hspace{4cm}\leq \frac{m!}{2}\left( \frac{\sqrt{17}}{c_1\wedge 1}\sqrt{\delta} \sigma_{\max}(1)\right)^2 \left( \frac{2\sigma_{\max}(1)}{c_1\wedge 1} \right)^{m-2}.
\end{align*}

Then, for any $\theta\in \mathcal S^p_s$, Bernstein's inequality gives for any $t'>0$ that
\begin{align*}
\P\left( \left|\theta^\top \left( \mathrm{diag}(\Sigma^{(\delta)}_n - \Sigma^{(\delta)})\right)\theta\right|> \frac{\sqrt{34}\sigma_{\max}(1)}{c_1\wedge 1}\sqrt{\frac{\delta t'}{n}} + \frac{2\sigma_{\max}(1)}{c_1\wedge 1}\frac{t'}{n} \right)  \leq 2e^{-t'}
\end{align*}
Next, a similar union bound argument as we used above for $\mathbf Z_n^{(2)}(s)$ gives
\begin{align}\label{bernstein-Z1}
\P\left( \mathbf Z_n^{(1)}(s) > \mathbf \zeta_n^{(1)}(s,t)  \right)\leq 2e^{-t},
\end{align}
with
\begin{align*}
\zeta_n^{(1)}(s,t)= \frac{\sigma_{\max}(1)}{c_1\wedge 1}\left(\sqrt{\frac{34\delta(t + s\log( 9) +  s\log\left(\frac{ep}{s}\right))}{n}}+\frac{2(t + s\log( 9) +  s\log\left(\frac{ep}{s}\right))}{n}\right).
\end{align*}
Next, easy computations give $\frac{1}{\delta}\zeta_n^{(1)}(s,t) + \frac{1}{\delta^2}\zeta_n^{(2)}(s,t)\leq\bar{\zeta}_n(s,t) $ where
\begin{align*}
\bar{\zeta}_n(s,t) &=\frac{\sigma_{\max}(s)}{c_1\wedge 1} \left( (16\sqrt{2}+ \sqrt{34}) \sqrt{\frac{t + s\log( 9) +  s\log\left(\frac{e p}{s}\right)}{\delta^{2} n}} + 10\frac{t + s\log( 9) +  s\log\left(\frac{ep}{s}\right)}{\delta^{2} n}  \right)
\end{align*}
Combining (\ref{proc-Z}), (\ref{bernstein-Z1}) and (\ref{bernstein-Z2}) with a union bound argument, we get, for any $s=1,\ldots,p$, that
\begin{align*}
\P\left( \mathbf Z_n(s) > \bar{\zeta}_n(s,t)  \right)\leq 4e^{-t}.
\end{align*}
Finally, using again a union bound argument, we get from the previous display that
\begin{align*}
\P\left( \bigcap_{s=1}^p \left\lbrace \mathbf Z_n(s) > \bar{\zeta}_n(s,t) \right\rbrace \right)\leq 4 p e^{-t}.
\end{align*}
Replacing $t$ by $t+\log(ep)$ and up to a rescaling of the constants, we get that
\begin{align*}
\P\left( \bigcap_{s=1}^p \left\lbrace \mathbf Z_n(s) > C\frac{\sigma_{\max}(s)}{ c_1\wedge 1}\max\left(\sqrt{\frac{t + s\log\left(ep/s\right)}{\delta^2 n}},\frac{t + s\log\left(ep/s\right)}{\delta^2 n}\right)  \right\rbrace \right)\leq e^{-t},
\end{align*}
for some numerical constant $C>0$.

\end{proof}

\subsection{Proof of Theorem \ref{main1}}

We will use the following lemma in order to prove our results
\begin{lemma}\label{curvature}
Let $\theta\in \mathcal S^p$. Let $\Sigma\in \R^{p\times p}$ be a symmetric positive semi-definite matrix with largest eigenvalue $\sigma_1$ of multiplicity $1$ and second largest eigenvalue $\sigma_2$. Then, for any $\theta\in \mathcal S^p$, we have
\begin{equation*}
\frac{1}{2}(\sigma_1 -\sigma_2)\|\theta\theta^\top - \theta_1\theta_1^\top \|_2^2 \leq \langle \Sigma,\theta_1\theta_1^\top -  \theta\theta^\top  \rangle.
\end{equation*}
\end{lemma}
See Lemma 3.2.1 in \cite{VL12} for a proof of this result.

\begin{proof}
We have by definition of $\hat\theta_1$ and in view of Lemma \ref{curvature} that
\begin{align*}
\frac{\sigma_1 - \sigma_2}{2}\| \hat\theta_1\hat\theta_1^\top - \theta_1\theta_1^\top  \|_2^2 &\leq \left\langle \Sigma, \theta_1\theta_1^\top - \hat\theta_1\hat\theta_1^\top \right\rangle\\
&\leq \left\langle \Sigma - \tilde \Sigma_n, \theta_1\theta_1^\top - \hat\theta_1\hat\theta_1^\top \right\rangle +  \left\langle \tilde \Sigma_n ,  \theta_1\theta_1^\top - \hat\theta_1\hat\theta_1^\top \right\rangle\\
&\leq \left\langle \Sigma - \tilde \Sigma_n, \theta_1\theta_1^\top -\hat\theta_1\hat\theta_1^\top \right\rangle  + \left[ \theta_1^\top \tilde \Sigma_n \theta_1 -\lambda |\theta_1|_0 \right]\\
&\hspace{3cm} -\left[ \hat\theta_1^\top  \tilde\Sigma_n\hat\theta_1- \lambda |\hat\theta_1|_0\right] +\lambda |\theta_1|_0-\lambda |\hat\theta_1|_0\\
&\leq \left\langle \Sigma - \tilde \Sigma_n, \theta_1\theta_1^\top -\hat\theta_1\hat\theta_1^\top \right\rangle  +\lambda |\theta_1|_0-\lambda |\hat\theta_1|_0\\
&\leq \|\Pi_{\hat J \cup J_1} (\Sigma - \tilde \Sigma_n)\Pi_{\hat J \cup J_1}\|_{\infty} \sqrt{2}\| \theta_1\theta_1^\top -\hat\theta_1\hat\theta_1^\top \|_2\\
&\hspace{6cm}  +\lambda |\theta_1|_0-\lambda |\hat\theta_1|_0,
\end{align*}
where $\Pi_{\hat J \cup J_1}$ is the orthogonal projection onto $\mathrm{l.s.}(e_j,\, j\in \hat J \cup J_1)$, $\hat J = J(\hat\theta_1)$ and $J_1= J(\theta_1)$.

Thus we get
\begin{align*}
\| \hat\theta_1\hat\theta_1^\top - \theta_1\theta_1^\top  \|_2^2 &\leq  \frac{2\sqrt{2}}{\sigma_1-\sigma_2}\|\Pi_{\hat J \cup J_1} (\Sigma - \tilde \Sigma_n)\Pi_{\hat J \cup J_1}\|_{\infty} \| \theta_1\theta_1^\top -\hat\theta_1\hat\theta_1^\top \|_2\\
&\hspace{5cm}   +\frac{2}{\sigma_1-\sigma_2}\lambda \left(|\theta_1|_0- |\hat\theta_1|_0\right).
\end{align*}
Set $A = \| \hat\theta_1\hat\theta_1^\top - \theta_1\theta_1^\top
\|_2 $, $\beta = \frac{2\sqrt{2}}{\sigma_1-\sigma_2}\|\Pi_{\hat J
\cup J_1} (\Sigma - \tilde \Sigma_n)\Pi_{\hat J \cup J_1}\|_{\infty}
$ and $\gamma =  \frac{2}{\sigma_1-\sigma_2}\lambda
\left(|\theta_1|_0- |\hat\theta_1|_0\right)$. The above display
becomes $$A^2 - \beta A - \gamma\leq 0.$$ Next, basic computations
on second order polynoms yield the following necessary condition on
$A$
$$A\leq \frac{\beta + \sqrt{\beta^2 + 4 \gamma}}{2} \leq \sqrt{\frac{2\beta^2 + 4\gamma}{2}} = \sqrt{\beta^2 + 2 \gamma}.$$
where we have used concavity of $x\rightarrow \sqrt{x}$.

Set $\hat s_1 = |\hat \theta_1|_0$ and $s_1 = |\theta_1|_0$.
Next, we have
\begin{align}\label{proof-thm-interm1}
\beta^2 + 2\gamma &\leq \frac{16}{(\sigma_1 - \sigma_2)^2} \left( \|\Pi_{\hat J}\left( \tilde \Sigma_n - \Sigma \right)\Pi_{\hat J}\|_\infty^2 - \frac{\sigma_1 - \sigma_2}{8}\lambda \hat s_1 \right)\notag\\
&\hspace{3cm} + \frac{16}{(\sigma_1 - \sigma_2)^2} \|\Pi_{ J_1}\left( \tilde \Sigma_n - \Sigma \right)\Pi_{ J_1}\|_\infty^2 +\frac{2\lambda}{\sigma_1-\sigma_2} s_1\notag\\
&\leq \frac{16}{(\sigma_1 - \sigma_2)^2} \max_{1\leq s \leq \bar s}\left\lbrace [\mathbf{Z}_n(s)]^2-\frac{(\sigma_1 - \sigma_2)}{8}\lambda s \right\rbrace\notag\\
&\hspace{5cm} + \frac{16}{(\sigma_1 - \sigma_2)^2} [\mathbf{Z}_n(s_1)]^2 +\frac{2\lambda}{\sigma_1-\sigma_2} s_1.
\end{align}

Next, we have in view of Proposition \ref{Prop-Zn} and under the condition $\bar s \log (ep/\bar s) \leq \delta^2 n$, with probability at least $1-\frac{1}{p}$ that
\begin{align*}
[\mathbf{Z}_n(s)]^2 \leq c^2 \frac{\sigma_{\max}^2(s)}{c_1^2 \wedge 1}\frac{(s+1)\log(ep/s)}{\delta^2 n}   ,\quad \forall 1\leq s \leq \bar s.
\end{align*}
Thus, we get with probability at least $1-\frac{1}{p}$ that
$$
[\mathbf{Z}_n(s_1)]^2 \leq c^2 \frac{\sigma_{\max}^2(s_1)}{c_1^2 \wedge 1}\frac{(s_1+1)\log(ep/s_1)}{\delta^2 n}.
$$
and
$$
\max_{1\leq s \leq \bar s}\left( [\mathbf{Z}_n(s)]^2-\frac{\sigma_1 - \sigma_2}{8}\lambda s \right) \leq 0,
$$
if we take
$$
\lambda = C \frac{\sigma_1^2}{\sigma_1 - \sigma_2}\frac{\log ep}{\delta^2 n},
$$
where $C>0$ is a large enough numerical constant.

Combining the last three displays with (\ref{proof-thm-interm1}), we get with probability at least $1-\frac{1}{p}$ that
\begin{align}
\beta^2 + 2\gamma &\leq C' \frac{\sigma_1^2}{(\sigma_1 - \sigma_2)^2} s_1 \frac{\log(ep)}{\delta^2 n},
\end{align}
where $C'>0$ is a numerical constant.
\end{proof}

\subsection{Proof of Lemma \ref{lem-datadriven}}

\begin{proof}
A standard matrix perturbation argument gives $|\hat\sigma_j - \sigma_j| \leq \| \tilde \Sigma_n - \Sigma \|_\infty$, $\forall 1\leq j \leq p$. Consequently, we get
\begin{align*}
\sigma_1 - \| \tilde \Sigma_n - \Sigma \|_\infty \leq  \hat \sigma_1  \leq  \sigma_1 + \| \tilde \Sigma_n - \Sigma \|_\infty,
\end{align*}
\begin{align*}
\hat \sigma_1 - \hat \sigma_2 &= \hat\sigma_1 -\sigma_1 + \sigma_1- \sigma_2 + \sigma_2 - \hat\sigma_2\\
&\geq \sigma_1- \sigma_2 - (|\hat\sigma_1 -\sigma_1| + |\hat\sigma_2 - \sigma_2|)\\
&\geq \sigma_1- \sigma_2 - 2\|\tilde \Sigma_n -\Sigma\|_\infty
\end{align*}
and similarly
\begin{align*}
\hat \sigma_1 - \hat \sigma_2 &\leq \sigma_1- \sigma_2 + 2\|\tilde \Sigma_n -\Sigma\|_\infty.
\end{align*}
Combining now (\ref{spectrum-bound}) with (\ref{nlowrank2}) with a sufficiently large constant $c>0$, we get with probability at least $1-\frac{1}{p}$ that
\begin{align*}
\frac{1}{2}\sigma_1 \leq \hat\sigma_1 \leq 2 \sigma_1, \quad \frac{1}{2}(\sigma_1-\sigma_2) \leq \hat\sigma_1 -\hat\sigma_2 \leq 2 (\sigma_1-\sigma_1),
\end{align*}
and
\begin{align*}
\frac{\sigma_1^2}{8(\sigma_1 - \sigma_2)} \leq \frac{\hat\sigma_1^2}{\hat\sigma_1 - \hat\sigma_2} \leq \frac{8\sigma_1^2}{\sigma_1 - \sigma_2}.
\end{align*}
The conclusion follows immediately.

\end{proof}

\subsection{Proof of Theorem \ref{thlowerbound}}

This proof uses standard tools of the minimax theory (cf. for
instance \cite{tsy_09}). The proof is more technical in the missing
observation case ($\delta<1$) in order the get the sharp dependence
$\delta^{-2}$ factor. In order to improve readability, we will
decompose the proof into several technical facts and proceed first
with the main arguments. Then, we give the proofs for the technical
facts.

\begin{proof}
We consider the following class $\tilde{ \mathcal{C}}$ of $p\times
p$ covariance matrices
\begin{align}\label{classsigma}
\tilde{\mathcal{C}} = \left\lbrace  \Sigma_\theta =
\Sigma(\theta,\sigma_1,\sigma_2)  = \sigma_1 \theta \theta^\top +
\sigma_2 (I_p -\theta\theta^\top),\,\forall \theta\in \mathcal
S^p\,:\,|\theta|_0\leq s_1,\, \forall \sigma_1\geq (1+\eta)\sigma_2
> 0\right\rbrace,
\end{align}
where $I_p$ is the $p\times p$ identity matrix and $\eta>0$ is some
absolute constant.

Note that the set $\tilde{\mathcal {C}}$ contains only full rank
matrices with the same determinant and whose first principal
component $\theta$ is $s_1$-sparse. Note also that $\tilde
{\mathcal{ C}} \subset \mathcal C $. Indeed, it is easy to see that
$\sigma_1$ is the largest eigenvalue of $\Sigma$ with multiplicity
$1$ and associated eigenvector $\theta$ with less than $s_1$ nonzero
components,
 $ \|I_p -\theta\theta^\top\|_\infty = 1$ and $(I_p -\theta\theta^\top)\theta = 0$.

Next, we define $\omega_0 = \left(1,1,0,\cdots,0\right)\in \{0,1\}^p$ and
$$
\Omega = \left\{\omega=(\omega^{(1)},\cdots,\omega^{(p)})\in \{0,1\}^p\,:\, \omega^{(1)} = \omega^{(2)}=1,\; |\omega|_0 =  s_1 \right\rbrace\cup \{  \omega_0 \}.
$$
A Varshamov-Gilbert's type bound (see for instance Lemma 4.10 in \cite{massart}) guarantees the existence of a subset $\mathcal N \subset \Omega$ with cardinality
$\log(\mathrm{Card}(\tilde{ \mathcal N })) \geq C_1 (s_1-2) \log(e(p-2)/(s_1-2))$ containing $\omega_0$ such that, for any two distinct elements $\omega$ and $\omega'$ of $\mathcal N$, we have
\begin{align*}
|\omega - \omega'|_0 \geq  \frac{s_1}{8}
\end{align*}
where $C_1>0$ is an absolute constant.

Set $\epsilon = a\sqrt{\frac{  \bar\sigma^2  s_1
\log(ep/s_1)}{\delta^2 n}}$ for some numerical constant $a\in
(0,1/\sqrt{2})$. Note that we have $\epsilon < 1/2$ under Condition
(\ref{nthreshold}). Consider now the following set of normalized
vectors
\begin{align}
\Theta &= \left\lbrace \theta(\omega) = \left(\sqrt{\frac{1-\epsilon^2}{2}},\sqrt{\frac{1-\epsilon^2}{2}},\frac{\omega^{(3)}\epsilon}{\sqrt{s_1-2}},\cdots, \frac{\omega^{(p)}\epsilon}{\sqrt{s_1-2}}\right)^\top\;:\; \omega\in \mathcal N\setminus\{\omega_0\}\right\rbrace\notag\\
&\hspace{2cm}\cup \{\theta_0 =\frac{1}{\sqrt{2}} \omega_0^\top\}.
\end{align}
Note that $|\Theta| = |\mathcal N|$ and $|\theta|_0 \leq s_1$ for any $\theta \in \Theta$.
\medskip

\begin{lemma}\label{distance projection}
For any $a>0$ and any distinct $\theta_1, \theta_2\in \Theta$, we have
\begin{align}\label{lower_1}
\| \theta_1\theta_1^\top - \theta_2\theta_2^\top \|_2^2 \geq \frac{a^2}{8}\bar\sigma^2 \frac{s_1\log(ep/s_1)}{\delta^2 n}.
\end{align}
\end{lemma}

\medskip
Clearly, for any $\theta\in \Theta$, we have $\Sigma_\theta \in \tilde{\mathcal {C}}$. We introduce now the class
$$\mathcal C(\Theta) = \left\{ \Sigma_\theta\in \tilde{\mathcal C} \;: \; \theta\in \Theta \right\}.$$
Denote by $\mathbb P_{\Sigma}$ the
distribution of $(Y_1,\cdots,Y_n)$. For any $\theta,\theta'\in \mathcal S^p$, the Kullback-Leibler
divergences $K\big(\P_{{
\Sigma_{\theta'}}},\P_{\Sigma_{\theta}}\big)$ between
$\P_{\Sigma_{\theta'}}$ and $\P_{\Sigma_{\theta}}$ is defined by
$$
K\left( \P_{\Sigma_{\theta'}},\P_{\Sigma_{\theta}} \right) = \mathbb E_{\Sigma_{\theta'}}\log \left(\frac{d\P_{\Sigma_{\theta'}}}{d\P_{\Sigma_{\theta}}}\right).
$$

We have the following result
\begin{lemma}\label{kullback-lem}
Let $X_1,\ldots,X_n\in\R^p$ be i.i.d.
$N\left(0,\Sigma\right)$ with $\Sigma =\Sigma_\theta\in \mathcal C(\Theta)$. Assume that $\frac{\sigma_1}{\sigma_2}\geq 1+\eta$ for some absolute $\eta>0$. Taking $a>0$ sufficiently small, we have for any $\theta' \in \mathcal S^p$, that
\begin{align*}
K\big(\P_{\Sigma_{\theta'}},\P_{\Sigma_{\theta_0}}\big)&\leq \frac{a^2}{2}s_1\log\left(\frac{ep}{s_1}\right).
\end{align*}
\end{lemma}

Thus, we have that
\begin{equation}\label{eq: condition C}
\frac{1}{\mathrm{Card}(\Theta)-1} \sum_{
\theta\in\Theta\setminus\{\theta_0\}}K(\P_{\Sigma_\theta}\P_{\Sigma_{\theta_0}})\
\leq\ \alpha \log \big(\mathrm{Card}(\Theta)-1\big)
\end{equation}
is satisfied for any $\alpha>0$ if $a>0$ is chosen as a sufficiently
small numerical constant depending on $\alpha$. In view of
(\ref{lower_1}) and (\ref{eq: condition C}), (\ref{eq:lower1}) now
follows by application of Theorem 2.5 in \cite{tsy_09}.
\end{proof}

\subsection{Proof of Lemma \ref{distance projection}}

\begin{proof}
For any distinct $\theta_1, \theta_2\in \Theta$, we have
\begin{align*}
|\theta_1 - \theta_2|_2^2 &\geq \frac{1}{8}\epsilon^2 =
\frac{a^2}{8}\bar\sigma^2 \frac{s_1\log(ep/s_1)}{\delta^2 n}.
\end{align*}

Next, we need to compare $\| \theta_1\theta_1^\top - \theta_2\theta_2^\top \|_2$ to $|\theta_1 - \theta_2|_2$.

For any $\theta_1,\theta_2\in \Theta$, we have
\begin{align}
\|\theta_1\theta_1^\top - \theta_2\theta_2^\top \|_2^2 &= 2 -2 (\theta_1^\top \theta_2)^2\notag\\
&= |\theta_1|_2^2 + |\theta_2|_2^2 -2 (\theta_1^\top \theta_2)^2 \notag\\
&= |\theta_1 -\theta_2|_2^2 +2[(\theta_1^\top \theta_2) - (\theta_1^\top \theta_2)^2]. \notag
\end{align}
We immediately get from the previous display that $\|
\theta_1\theta_1^\top - \theta_2\theta_2^\top \|_2 \geq |\theta_1 -
\theta_2|_2$ for any $\theta_1,\theta_2\in \Theta$.

\end{proof}

\subsection{Proof of Lemma \ref{kullback-lem}}

Recall that $X_1,\ldots,X_n\in\R^p$ are i.i.d.
$N\left(0,\Sigma\right)$ with $\Sigma =\Sigma_\theta\in \mathcal C(\Theta)$. For any $1\leq i \leq n$, set  $\delta_i = \left(\delta_{i,1},\cdots,\delta_{i,p}\right)^\top\in\R^p$. We note that $\delta_1,\ldots,\delta_n$ are random vectors in $\R^p$
with i.i.d. entries $\delta_{i,j} \sim B(\delta)$ and independent from $(X_1,\cdots,X_n)$. Recall that the
observations $Y_1,\ldots,Y_n$ satisfies $Y^{(j)}_i =
\delta_{i,j}X_{i}^{(j)}$. Denote by $\mathbb P_{\Sigma}$ the
distribution of $(Y_1,\cdots,Y_n)$ and by $\mathbb
P_{\Sigma}^{(\delta)}$ the conditional distribution of
$(Y_1,\cdots,Y_n)$ given $(\delta_{1},\cdots,\delta_n)$. Next, we
note that for any $1\leq i \leq n$ the conditional random variables
$Y_i\mid (\delta_{1},\cdots,\delta_{n})$ are independent Gaussian vectors
$N(0,\Sigma_\theta^{(\delta_i)})$, where
\begin{equation*}
(\Sigma_\theta^{(\delta_i)})_{j,k} =
\begin{cases}
\delta_{i,j}\delta_{i,k}\Sigma_{j,k}&\text{if $j\neq k$},\\
\delta_{i,j}\Sigma_{j,j}&\text{otherwise}.
\end{cases}
\end{equation*}

Thus, we have $\P_{\Sigma_{\theta}}^{(\delta)} = \otimes_{i=1}^n
\P_{\Sigma_{\theta}^{(\delta_i)}}$. Denote respectively by $\P_\delta$ and $\E_\delta$ the
probability distribution  of $(\delta_1,\cdots,\delta_n)$ and the
associated expectation. We also denote by $\E_{\Sigma_{\theta}}$
and $\E^{(\delta)}_{\Sigma_{\theta}}$ the expectation and conditional
expectation associated respectively with $\P_{\Sigma_{\theta}}$ and
$\P_{\Sigma_{\theta}}^{(\delta)}$.

Next, for any $\theta,\theta'\in \mathcal S^p$, the Kullback-Leibler
divergences $K\big(\P_{{
\Sigma_{\theta'}}},\P_{\Sigma_{\theta}}\big)$ between
$\P_{\Sigma_{\theta'}}$ and $\P_{\Sigma_{\theta}}$ satisfies
\begin{align}\label{proof-theo2-interm2}
K\left( \P_{\Sigma_{\theta'}},\P_{\Sigma_{\theta}} \right)& =
\mathbb E_{\Sigma_{\theta'}}\log
\left(\frac{d\P_{\Sigma_{\theta'}}}{d\P_{\Sigma_{\theta}}}\right) =
\mathbb E_{\Sigma_{\theta'}}\log \left(\frac{d(\P_{\delta}\otimes
\P_{\Sigma_{\theta'}}^{(\delta)})}{d(\P_{\delta}\otimes
\P_{\Sigma_{\theta}}^{(\delta)})}\right)
\notag\\
&=\E_{\delta}\E_{\Sigma_{\theta'}}^{(\delta)}\log \left(\frac{d
\P_{\Sigma_{\theta'}}^{(\delta)}}{d
\P_{\Sigma_{\theta}}^{(\delta)}}\right)= \E_{\delta}K\left( \P_{
\Sigma_{\theta'}}^{(\delta)},\P_{\Sigma_{\theta}}^{(\delta)}
\right)\notag\\
&= \sum_{i=1}^n \E_{\delta_i} K\left( \P_{
\Sigma_{\theta'}^{(\delta_i)}},\P_{\Sigma_{\theta}^{(\delta_i)}}
\right).
\end{align}

Set $\theta_{\delta_i} = (\delta_{i,1}\theta^{(1)},\cdots,\delta_{i,p}\theta^{(p)})^\top$. In view of (\ref{classsigma}), we have
\begin{align}\label{proof-theo2-interm1}
\Sigma_\theta^{(\delta_i)} =
\left[(\sigma_1-\sigma_2)|\theta_{\delta_i}|_2^2 +
\sigma_2\right]\Pi_{\theta,\delta_i} + \sigma_2 \left(
I_p^{(\delta_i)} - \Pi_{\theta,\delta_i} \right),
\end{align}
and $\Pi_{\theta,\delta_i}$ is the orthogonal projection onto $l.s.(\theta_{\delta_i})$
 (Note indeed that we have in general $|\theta_{\delta_i}|_2\leq 1$, therefore $\Pi_{\theta,\delta_i} = |\theta_{\delta_i}|_2^{-2} \theta_{\delta_i}\theta_{\delta_i}^\top$). For any $\theta\in \Theta$, we set $\sigma_1(\theta) = (\sigma_1-\sigma_2)|\theta_{\delta_i}|_2^2
+ \sigma_2$.

\medskip
\begin{description}
 \item \textbf{Fact 1:} For any $1\leq i \leq
n$, any $\theta,\theta' \in \mathcal S^p$ and any realization of
$\delta_i\in \{0,1\}^p$, we have
\begin{align*}
K\big(\P_{{\Sigma_{\theta'}^{(\delta_i)}}},\P_{\Sigma_{\theta}^{(\delta_i)}}\big) &= \frac{1}{2}\left(
\frac{\sigma_2}{\sigma_1(\theta)}+\frac{\sigma_1(\theta')}{\sigma_2}-2\right) + \frac{1}{2}\log\left(
\frac{\sigma_1(\theta)}{\sigma_1(\theta')}\right)\\
&\hspace{0.5cm}+\frac{1}{2}\mathrm{tr}\left(\Pi_{\theta,\delta_i}\Pi_{\theta',\delta_i}\right)\left[
\frac{\sigma_1(\theta')}{\sigma_1(\theta)}+1-\frac{\sigma_2}{\sigma_1(\theta)}-\frac{\sigma_1(\theta')}{\sigma_2}\right].
\end{align*}
\end{description}
\medskip

We apply Fact 1 with $\theta = \theta_0 = \frac{1}{\sqrt{2}}\omega_0^\top$ and take the expectation w.r.t. $\delta_i$ for any $i=1,\ldots,n$. Thus, we get the following.
\medskip
\begin{description}
 \item \textbf{Fact 2:} Assume that $\frac{\sigma_1}{\sigma_2}\geq 1+\eta$ for some absolute $\eta>0$. Taking $a>0$ sufficiently small (that can depend only on $\eta$), we have for any $i=1,\ldots,n$, any $\theta' \in \mathcal S^p$,
 that
\begin{align*}
\E_{\delta_i}\left[K\big(\P_{{\Sigma_{\theta'}^{(\delta_i)}}},\P_{\Sigma_{\theta_0}^{(\delta_i)}}\big)
\right]&\leq \frac{\delta^2}{2\bar \sigma^2}\epsilon^2.
\end{align*}
\end{description}
\medskip

We immediately get from Fact 2 for any $\theta' \in \Theta$ that
\begin{align*}
K\big(\P_{{\Sigma_{\theta'}}},\P_{\Sigma_{\theta_0}}\big)=\sum_{i=1}^n\E_{\delta_i}\left[K\big(\P_{{\Sigma_{\theta'}^{(\delta_i)}}},\P_{\Sigma_{\theta_0}^{(\delta_i)}}\big)
\right]& \leq \frac{\delta^2 n}{2\bar \sigma^2}\epsilon^2 =
\frac{a^2}{2} s_1 \log(ep/s_1).
\end{align*}

\subsection{Proof of Fact 1}

\begin{proof}
In view of (\ref{proof-theo2-interm1}), we get for any $1\leq i \leq
n$, any $\theta,\theta' \in \mathcal S^p$ and any realization
$\delta_i\in \{0,1\}^p$ that $\P_{\Sigma_{\theta}^{(\delta_i)}}\ll
\P_{\Sigma_{\theta'}^{(\delta_i)}}$ and hence
$K\big(\P_{{\Sigma_{\theta'}^{(\delta_i)}}},\P_{\Sigma_{\theta}^{(\delta_i)}}\big)<\infty$.

Define $J_i = \left\lbrace  j\,:\, \delta_{i,j}=1,\, 1\leq j \leq r
\right\rbrace$ and $d_i = |J_i|$. Define the mapping
$P_i:\R^p\rightarrow \R^{d_i}$ as follows $P_i(x)= x(J_i)$ where for
any $x=(x^{(1)},\cdots,x^{(p)})^\top\in\R^p$, $x(J_i)\in \R^{d_i}$
is obtained by keeping only the components $x^{(k)}$ with their
index $k\in J_i$. We denote  by $P_i^{*}\,: \R^{d_i} \rightarrow
\R^{p}$ the right inverse application of $P_i$. We note that
\begin{align*}
P_i \Sigma_{\theta}^{(\delta_i)}P_i^{*} = \sigma_1(\theta)\Pi_{\theta(J_i),\delta_i} + \sigma_2 \left[ I_{d_i} - \Pi_{\theta(J_i),\delta_i} \right],
\end{align*}
where $\Pi_{\theta(J_i),\delta_i}$ denotes the orthogonal
projection onto the subspace $\mathrm{l.s.}(\theta_{\delta_i}(J_i))$
of $\R^{d_i}$. Note also that $P_i
\Sigma_{\theta}^{(\delta_i)}P_i^{*}$ admits an inverse for any
$\theta\in \mathcal S^p$ provided that $\delta_i$ is not the null
vector in $\R^p$ and we have
$$
(P_i \Sigma_{\theta}^{(\delta_i)}P_i^{*})^{-1} = \frac{1}{\sigma_{1}(\theta)}\Pi_{\theta(J_i),\delta_i} + \frac{1}{\sigma_2} \left[ I_{d_i} - \Pi_{\theta(J_i),\delta_i} \right].
$$

Thus, we get for any $\theta,\theta'\in \mathcal S^p$ that
\begin{align*}
&\hspace{1cm} K\big(\P_{{\Sigma_{\theta'}^{(\delta_i)}}},\P_{\Sigma_{\theta}^{(\delta_i)}}\big)= K\big(\P_{P_i \Sigma_{\theta'}^{(\delta_i)}P_i^{*}},\P_{P_i(\Sigma_{\theta}^{(\delta_i)})P_i^{*}}\big)\\
&\hspace{0.5cm}= \frac{1}{2}\mathrm{tr}\left((P_i \Sigma_{\theta}^{(\delta_i)}P_i^{*})^{-1}  P_i(\Sigma_{\theta'}^{(\delta_i)})P_i^{*}\right) + \frac{1}{2}\log\left( \frac{\mathrm{det}\left(P_i\Sigma_{\theta}^{(\delta_i)}P_i^*\right)}{\mathrm{det}\left(P_i\Sigma_{\theta'}^{(\delta_i)}P_i^*\right)}\right) -\frac{d_i}{2}\\
&\hspace{0.5cm}= \frac{1}{2}\mathrm{tr}\left(\left[\frac{1}{\sigma_{1}(\theta)}\Pi_{\theta(J_i),\delta_i} + \frac{1}{\sigma_2} \left( I_{d_i} - \Pi_{\theta(J_i),\delta_i} \right)\right] \left[\sigma_1(\theta')\Pi_{\theta'(J_i),\delta_i} + \sigma_2 \left(I_{d_i} - \Pi_{\theta'(J_i),\delta_i} \right)\right]\right)\\
&\hspace{4cm}  + \frac{1}{2}\log\left( \frac{\mathrm{det}\left(P_i\Sigma_{\theta}^{(\delta_i)}P_i^*\right)}{\mathrm{det}\left(P_i\Sigma_{\theta'}^{(\delta_i)}P_i^*\right)}\right) -\frac{d_i}{2}\\
&\hspace{0.5cm}= \frac{1}{2}\left( \frac{\sigma_2}{\sigma_1(\theta)} + \frac{\sigma_1(\theta')}{\sigma_2} -2\right) + \frac{1}{2}\mathrm{tr}\left( \Pi_{\theta(J_i),\delta_i}\Pi_{\theta'(J_i),\delta_i}\right)\left[ \frac{\sigma_1(\theta')}{\sigma_1(\theta)}+  1 - \frac{\sigma_2}{\sigma_1(\theta)} - \frac{\sigma_1(\theta')}{\sigma_2} \right]\\
&\hspace{4cm}  + \frac{1}{2}\log\left( \frac{\mathrm{det}\left(P_i\Sigma_{\theta}^{(\delta_i)}P_i^*\right)}{\mathrm{det}\left(P_i\Sigma_{\theta'}^{(\delta_i)}P_i^*\right)}\right) \\
&\hspace{0.5cm}= \frac{1}{2}\left(
\frac{\sigma_2}{\sigma_1(\theta)}+\frac{\sigma_1(\theta')}{\sigma_2}-2\right)+\frac{1}{2}\mathrm{tr}\left(\Pi_{\theta,\delta_i}\Pi_{\theta',\delta_i}\right)\left[
\frac{\sigma_1(\theta')}{\sigma_1(\theta)}+1-\frac{\sigma_2}{\sigma_1(\theta)}-\frac{\sigma_1(\theta')}{\sigma_2}\right]\\
&\hspace{4cm}  + \frac{1}{2}\log\left(
\frac{\sigma_1(\theta)}{\sigma_1(\theta')}\right),
\end{align*}
where we have used $\sigma_1(\theta)$ and $\sigma_2$ are eigenvalues of $P_i\Sigma_{\theta}^{(\delta_i)}P_i^*$ with respective multiplicity $1$ and  $d_1-1$ for any $\theta\in \Theta$, and also that $\mathrm{tr}\left( \Pi_{\theta(J_i),\delta_i}\Pi_{\theta'(J_i),\delta_i}\right) = \mathrm{tr}\left( \Pi_{\theta,\delta_i}\Pi_{\theta',\delta_i}\right)$ for any $\theta,\theta'\in \Theta$.
\end{proof}

\subsection{Proof of Fact 2}
\begin{proof}

For any $i=1,\ldots,n$, we have that
$\sigma_1(\theta_0)  = (\sigma_1 - \sigma_2)
\frac{\delta_{i,1}+ \delta_{i,2}}{2} + \sigma_2$ and
$$\sigma_1(\theta') = (\sigma_1 - \sigma_2)
\left[\left(\frac{1-\epsilon^2}{2}\right)(\delta_{i,1} +
\delta_{i,2}) + \epsilon^2\sum_{j=3}^p \frac{\omega^{(j)}}{s_1-2}\delta_{i,j}
\right] + \sigma_2.
$$
Note that the random quantities in the last three displays depend on
$\delta_{i,1}, \delta_{i,2}$ only through the sum $Z_i
:=\delta_{i,1} + \delta_{i,2}\sim \mathrm{Bin}(2,\delta)$ and that $
Z_i$ is independent of $(\delta_{i,3},\cdots,\delta_{i,p})$.

Thus, we get
\begin{align*}
\E_{\delta_i}\left[\frac{\sigma_2}{2\sigma_1(\theta_0)}\right] & = \E_{\delta_i}\left[\frac{\sigma_2}{(\sigma_1-\sigma_2)(\delta_{i,1}+ \delta_{i,2}) + 2\sigma_2}\right]\\
&=\frac{\delta^2\sigma_2}{2\sigma_1} + \frac{2\delta(1-\delta)\sigma_2}{\sigma_1+ \sigma_2} + \frac{(1-\delta)^2}{2}.
\end{align*}
Similarly, we obtain
\begin{align*}
\E_{\delta_i}\left[\frac{\sigma_1(\theta')}{2\sigma_2}\right] &= \frac{\delta(\sigma_1-\sigma_2)|\theta'|_2^2+ \sigma_2}{ 2\sigma_2}=\frac{\delta\sigma_1}{2\sigma_2} + \frac{1-\delta}{2} .
\end{align*}
Combining the last two displays, we get
\begin{align}\label{annexinterm1}
\E_{\delta_i}\left[\frac{\sigma_2}{2\sigma_1(\theta_0)} +
\frac{\sigma_1(\theta')}{2\sigma_2} - 1\right] &=
\frac{\delta(\sigma_1-\sigma_2)^2(\sigma_1 +
\delta\sigma_2)}{2\sigma_1\sigma_2(\sigma_1+ \sigma_2)}.
\end{align}

We study now the following quantity
$$
\frac{1}{2}\mathrm{tr}\left(\Pi_{\theta_0,\delta_i}\Pi_{\theta',\delta_i}\right)\left[
\frac{\sigma_1(\theta')}{\sigma_1(\theta_0)}+1-\frac{\sigma_2}{\sigma_1(\theta_0)}-\frac{\sigma_1(\theta')}{\sigma_2}\right].
$$

we note first that
\begin{align*}
\Pi_{\theta_0,\delta_i} = \frac{1}{\delta_{i,1}+\delta_{i,2}}
\left(\begin{array}{cc|c} \delta_{i,1} &
\delta_{i,1}\delta_{i,2}&O\\
\delta_{i,1}\delta_{i,2} &
\delta_{i,2}&O\\
\hline O&O &O
\end{array}\right)
\end{align*}
and
\begin{align*}
\Pi_{\theta',\delta_i} =
\frac{1-\epsilon^2}{2|\theta'_{\delta_i}|_2^2}
\left(\begin{array}{cc|c} \delta_{i,1}&
\delta_{i,1}\delta_{i,2}&\ast\\
\delta_{i,1}\delta_{i,2} &\delta_{i,2}&\ast\\
\hline \ast&\ast &\ast
\end{array}\right)
\end{align*}
Thus, we get that
\begin{align*}
\frac{1}{2}\mathrm{tr}\left(\Pi_{\theta_0,\delta_i}
\Pi_{\theta',\delta_i}  \right) &=
\left(1-\epsilon^2\right)\frac{\delta_{i,1}^2 +
2\delta_{i,1}\delta_{i,2} +
\delta_{i,2}^2}{4|\theta'_{\delta_i}|_2^2(\delta_{i,1} +
\delta_{i,2})}\\
&=\left(1-\epsilon^2\right)\frac{\delta_{i,1}+\delta_{i,2}}{4|\theta'_{\delta_i}|_2^2}.
\end{align*}
Next, we set
\begin{align*}
\tilde{\sigma}(\theta_0,\theta')=\frac{\sigma_2 \sigma_1(\theta') +
\sigma_2\sigma_1(\theta_0) - \sigma_2^2 -
\sigma_1(\theta_0)\sigma_1(\theta')}{\sigma_2 \sigma_1(\theta_0)}.
\end{align*}

If $Z_i=1$, then $\sigma_1(\theta_0) = (\sigma_1 + \sigma_2)/2$ and
\begin{align*}
\tilde\sigma(\theta_0,\theta') &= \frac{2\sigma_2\sigma_1(\theta') +
\sigma_2(\sigma_1 + \sigma_2) -2\sigma_2^2 - (\sigma_1 +
\sigma_2)\sigma_1(\theta')}{\sigma_2(\sigma_1+\sigma_2)}\\
&=-\frac{(\sigma_1-\sigma_2)^2|\theta'_{\delta_i}|_2^2}{\sigma_2(\sigma_1+\sigma_2)}.
\end{align*}
If $Z_i=2$, then $\sigma_1(\theta_0) = \sigma_1$ and
\begin{align*}
\tilde\sigma(\theta_0,\theta') &=
-\frac{(\sigma_1-\sigma_2)^2|\theta'_{\delta_i}|_2^2}{\sigma_1\sigma_2}.
\end{align*}

We now freeze $(\delta_{i,3},\cdots,\delta_{i,p})$ and compute the
following expectation w.r.t $Z_i$
\begin{align*}
&\mathbb E_{Z_i}\left[  \frac{1}{2}\mathrm{tr}\left(
\Pi_{\theta_0,\delta_i} \Pi_{\theta',\delta_i}  \right)
\tilde{\sigma}(\theta_0,\theta') \right] =
-\delta^2\left(1-\epsilon^2\right)\frac{(\sigma_1-\sigma_2)^2}{2\sigma_1\sigma_2}\\
&\hspace{6cm}-\delta(1-\delta)\left(1-\epsilon^2\right)\frac{(\sigma_1-\sigma_2)^2}{2\sigma_2(\sigma_1+\sigma_2)}.
\end{align*}
We note that the above display does not depend on
$(\delta_{i,3},\cdots,\delta_{i,p})$. Thus, we get
\begin{align*}
&\mathbb E_{\delta}\left[  \frac{1}{2}\mathrm{tr}\left(
\Pi_{\theta_0,\delta_i} \Pi_{\theta',\delta_i}  \right)
\tilde{\sigma}(\theta_0,\theta') \right] =
-\delta^2\left(1-\epsilon^2\right)\frac{(\sigma_1-\sigma_2)^2}{2\sigma_1\sigma_2}\\
&\hspace{6cm}-\delta(1-\delta)\left(1-\epsilon^2\right)\frac{(\sigma_1-\sigma_2)^2}{2\sigma_2(\sigma_1+\sigma_2)}.
\end{align*}
Combining the above display with (\ref{annexinterm1}), we get
\begin{align*}
\Delta_1 &:=\frac{\delta(\sigma_1-\sigma_2)^2(\sigma_1 +
\delta\sigma_2)}{2\sigma_1\sigma_2(\sigma_1+ \sigma_2)}
-\delta^2\left(1-\epsilon^2\right)\frac{(\sigma_1-\sigma_2)^2}{2\sigma_1\sigma_2}\\
&\hspace{6cm}-\delta(1-\delta)\left(1-\epsilon^2\right)\frac{(\sigma_1-\sigma_2)^2}{2\sigma_2(\sigma_1+\sigma_2)}\\
&=-\frac{\delta}{2}\frac{(\sigma_1 - \sigma_2)^2}{\sigma_1
\sigma_2(\sigma_1 + \sigma_2)}\left(
(1-\delta)(1-2\epsilon^2)\sigma_1 - \delta \epsilon^2(\sigma_1 +
\sigma_2)\right).
\end{align*}

We study now the logarithm factor
$$
\Delta_2:=  \frac{1}{2}\E_\delta\log\left(
\frac{\sigma_1(\theta_0)}{\sigma_1(\theta')}\right).$$

Recall that $\sigma_1(\theta_0) = (\sigma_1 - \sigma_2)\frac{Z_i}{2}
+ \sigma_2$ with $Z_i = \delta_{i,1} + \delta_{i,2}\sim
\mathrm{Bin}(2,\delta)$ and
\begin{align*}
\sigma_1(\theta') &=(\sigma_1 -\sigma_2) |\theta_{\delta_i}|_2^2
+\sigma_2 = (\sigma_1 -\sigma_2) \left[ \frac{Z_i}{2} +
\frac{\epsilon^2}{s_1-2} \tilde Z_i \right] + \sigma_2.
\end{align*}
with $\tilde Z_i = \sum_{j=3}^p \omega^{(j)}\delta_{i,j} \sim
\mathrm{Bin}(s_1-2,\delta)$ and is independent of
$(\delta_{i,1},\delta_{i,2})$.

We now freeze $\tilde Z_i$ and take the expectation w.r.t. $Z_{i}$.
Thus, we get
\begin{align*}
\E_{Z_i}\left[
\frac{1}{2}\log\left(\frac{\sigma_1(\theta_0)}{\sigma_1(\theta')}\right)
\right] &= -\frac{(1-\delta)^2}{2}\log\left(\frac{(\sigma_1
-\sigma_2)\frac{\epsilon^2}{s_1-2} \tilde Z_i+ \sigma_2}{\sigma_2} \right) \\
&- \delta(1- \delta) \log\left(
\frac{(\sigma_1-\sigma_2)[(1-\epsilon^2) + \frac{2\epsilon^2}{s_1-2}
\tilde Z_i] +2\sigma_2}{\sigma_1 + \sigma_2} \right)\\ &-
\frac{\delta^2}{2} \log
\left(\frac{(\sigma_1-\sigma_2)[(1-\epsilon^2)+
\frac{\epsilon^2}{s_1-2} \tilde Z_i] +\sigma_2}{\sigma_1}\right).
\end{align*}
We study now the first term in the right-hand side of the above
display. We have
\begin{align*}
\frac{(\sigma_1 -\sigma_2)\frac{\epsilon^2}{s_1-2} \tilde Z_i+
\sigma_2}{\sigma_2}  &= 1-\epsilon^2 + \epsilon^2 \left[\left(
\frac{\sigma_1}{\sigma_2} -1  \right)\frac{\tilde Z_i}{s_1-2}
+1\right]\\
&=1-\epsilon^2 + \epsilon^2 \left[\sum_{j=3}^p \left[ \left(
\frac{\sigma_1}{\sigma_2} -1  \right)\delta_{i,j} +
1\right]\frac{\omega^{(j)}}{s_1-2}\right],
\end{align*}
since $\sum_{j=3}^p \omega^{(j)} = s_1-2$ by construction. Next, we
notice that $-\log$ is convex. Thus applying Jensen's inequality
twice gives that
\begin{align*}
\E_{\tilde Z_i} \left[-\log \left(\frac{(\sigma_1
-\sigma_2)\frac{\epsilon^2}{s_1-2} \tilde Z_i+ \sigma_2}{\sigma_2}
\right)\right] &\leq \epsilon^2\E_{\tilde Z_i} \left[-\log
\left(\sum_{j=3}^p \left[ \left( \frac{\sigma_1}{\sigma_2} -1
\right)\delta_{i,j} +
1\right]\frac{\omega^{(j)}}{s_1-2}\right)\right]\\
&\leq -\epsilon^2 \sum_{j=3}^p
\frac{\omega^{(j)}}{s_1-2}\E_{\delta_{i,j}}\log\left[\left(
\frac{\sigma_1}{\sigma_2} -1 \right)\delta_{i,j} + 1 \right]\\
&\leq  - \delta \epsilon^2
\log\left(\frac{\sigma_1}{\sigma_2}\right).
\end{align*}
We proceed similarly and obtain
\begin{align*}
\frac{(\sigma_1-\sigma_2)[(1-\epsilon^2) + \frac{2\epsilon^2}{s_1-2}
\tilde Z_i] +2\sigma_2}{\sigma_1 + \sigma_2} &=(1-\epsilon^2) +
\epsilon^2 \left(\frac{2\sigma_2 }{\sigma_1+ \sigma_2} +
\frac{2(\sigma_1-\sigma_2)\tilde Z_i}{(s_1-2)(\sigma_1 +
\sigma_2)}\right)\\
&=(1-\epsilon^2) + \epsilon^2 \sum_{j=3}^p
\frac{\omega^{(j)}}{s_1-2} \left(\frac{2\sigma_2
+\delta_{i,j}(\sigma_1-\sigma_2) }{\sigma_1+ \sigma_2}\right),
\end{align*}
and
\begin{align*}
&\E_{\tilde Z_i} \left[-\log
\left(\frac{(\sigma_1-\sigma_2)[(1-\epsilon^2) +
\frac{2\epsilon^2}{s_1-2} \tilde Z_i] +2\sigma_2}{\sigma_1 +
\sigma_2} \right)\right]\\
&\hspace{4cm}\leq -\epsilon^2 \sum_{j=3}^p
\frac{\omega^{(j)}}{s_1-2}\E_{\delta_{i,j}}\log\left[\frac{2\sigma_2
+\delta_{i,j}(\sigma_1-\sigma_2) }{\sigma_1+ \sigma_2} \right]\\
&\hspace{4cm}\leq  - (1-\delta) \epsilon^2 \log\left(\frac{2\sigma_2}{\sigma_1
+\sigma_2}\right).
\end{align*}
We obtain similarly
\begin{align*}
\frac{(\sigma_1-\sigma_2)[(1-\epsilon^2)+ \frac{\epsilon^2}{s_1-2}
\tilde Z_i] +\sigma_2}{\sigma_1} &=(1-\epsilon^2) + \epsilon^2
\left(\frac{\tilde Z_i(\sigma_1-\sigma_2)}{(s_1-2)\sigma_1}
+\frac{\sigma_2 }{\sigma_1}
\right)\\
&=(1-\epsilon^2) + \epsilon^2 \sum_{j=3}^p
\frac{\omega^{(j)}}{s_1-2} \left(\frac{\sigma_2
+\delta_{i,j}(\sigma_1-\sigma_2)}{\sigma_1}\right),
\end{align*}
and
\begin{align*}
&\E_{\tilde Z_i} \left[-\log
\left(\frac{(\sigma_1-\sigma_2)[(1-\epsilon^2)+
\frac{\epsilon^2}{s_1-2} \tilde Z_i] +\sigma_2}{\sigma_1}
\right)\right]\\
&\hspace{4cm}\leq -\epsilon^2 \sum_{j=3}^p
\frac{\omega^{(j)}}{s_1-2}\E_{\delta_{i,j}}\log\left[\frac{\sigma_2
+\delta_{i,j}(\sigma_1-\sigma_2)}{\sigma_1} \right]\\
&\hspace{4cm}\leq  - (1-\delta) \epsilon^2
\log\left(\frac{\sigma_2}{\sigma_1}\right).
\end{align*}
Thus, we get
\begin{align*}
\Delta_2 &= - \frac{(1-\delta)^2\delta}{2} \epsilon^2
\log\left(\frac{\sigma_1}{\sigma_2}\right)
-\delta(1-\delta)^2\epsilon^2 \log\left(\frac{2\sigma_2}{\sigma_1
+\sigma_2}\right) - \frac{(1-\delta)\delta^2}{2} \epsilon^2
\log\left(\frac{\sigma_2}{\sigma_1}\right)\\
&= - \frac{(1-\delta)\delta}{2} \epsilon^2 \left[
(1-\delta)\log\left(\frac{\sigma_1}{\sigma_2}\right) + 2
\log\left(\frac{2\sigma_2}{\sigma_1 +\sigma_2}\right) +\delta
\log\left(\frac{\sigma_2}{\sigma_1}\right) \right]\\
&= \frac{(1-\delta)\delta}{2} \epsilon^2 \left[
\log\left(\frac{(\sigma_1 + \sigma_2)^2}{2\sigma_1\sigma_2}\right)
-\log(2) + 2\delta \log\left(\frac{\sigma_1}{\sigma_2}\right)
\right]
\end{align*}
Set $\Delta:=\Delta_1+\Delta_2$. We have
\begin{align*}
\Delta &=-\frac{\delta}{2}\frac{(\sigma_1 - \sigma_2)^2}{\sigma_1
\sigma_2(\sigma_1 + \sigma_2)}\left(
(1-\delta)(1-2\epsilon^2)\sigma_1 - \delta \epsilon^2(\sigma_1 +
\sigma_2)\right)\\
&\hspace{4cm} + \frac{(1-\delta)\delta}{2} \epsilon^2 \left[
\log\left(\frac{(\sigma_1 + \sigma_2)^2}{4\sigma_1\sigma_2}\right)+
2\delta \log\left(\frac{\sigma_1}{\sigma_2}\right)
\right]\\
&\leq \frac{\delta^2}{2\bar \sigma^2}\epsilon^2 -\frac{\delta(1-
\delta)}{2}\left[ \frac{(\sigma_1 - \sigma_2)^2}{ \sigma_2(\sigma_1
+ \sigma_2)}\left((1-2\epsilon^2)\right) -\epsilon^2 \left[
\log\left(\frac{(\sigma_1 + \sigma_2)^2}{4\sigma_1\sigma_2}\right) +
2 \log\left(\frac{\sigma_1}{\sigma_2}\right) \right] \right]\\
&\leq\frac{\delta^2}{2\bar \sigma^2}\epsilon^2 -\frac{\delta(1-
\delta)}{2}\left[ \frac{(\sigma_1 - \sigma_2)^2}{ \sigma_2(\sigma_1
+ \sigma_2)}\left((1-2\epsilon^2)\right) -\epsilon^2
\log\left(\frac{\sigma_1(\sigma_1 + \sigma_2)^2}{4\sigma_2^3}\right)
\right].
\end{align*}
We now show that if the absolute constant $a>0$ (recall that
$\epsilon=a\bar\sigma\sqrt{\frac{s_1\log(ep/s_1)}{\delta^2 n}}$) is
taken sufficiently small, then we have
\begin{align*}
 \frac{(\sigma_1 - \sigma_2)^2}{ \sigma_2(\sigma_1
+ \sigma_2)}\left((1-2\epsilon^2)\right) -\epsilon^2
\log\left(\frac{\sigma_1(\sigma_1 +
\sigma_2)^2}{4\sigma_2^3}\right)\geq 0,\quad \forall \sigma_1
>(1+\eta)\sigma_2.
\end{align*}
We set $u = \sigma_1 -\sigma_2$ and $x=u/(2\sigma_2)$. Then, we have
\begin{align*}
& \frac{(\sigma_1 - \sigma_2)^2}{ \sigma_2(\sigma_1 +
\sigma_2)}\left((1-2\epsilon^2)\right) -\epsilon^2
\log\left(\frac{\sigma_1(\sigma_1 + \sigma_2)^2}{4\sigma_2^3}\right)\\
&\hspace{1cm}= \frac{u^2}{\sigma_2(2\sigma_2 +u)}(1-2\epsilon^2) - \epsilon^2
\log\left( \frac{(\sigma_2+u)(2\sigma_2+u)^2}{4\sigma_2^3} \right)\\
&\hspace{1cm}=
\frac{u^2}{2\sigma_2^2(1+u/(2\sigma_2))}(1-2\epsilon^2) -\epsilon^2
\log\left(\left[1+\frac{u}{\sigma_2}\right]\left[1+\frac{u}{2\sigma_2}\right]^2\right)\\
&\hspace{1cm}\geq
\frac{u^2}{2\sigma_2^2(1+u/(2\sigma_2)}(1-2\epsilon^2) -3\epsilon^2
\log\left(1+\frac{u}{\sigma_2}\right)\\
&\hspace{1cm}\geq  2(1-2\epsilon^2)\frac{x^2}{1+x} -3\epsilon^2
\log\left(1+2x\right)\geq 0,\quad \forall x\geq \frac{\eta}{2},
\end{align*}
provided that the numerical constant $a>0$ is taken sufficiently
small (and this choice can depend only on $\eta>0$).

\end{proof}

 \footnotesize{
\bibliographystyle{plain}
\bibliography{panel}
}

\end{document}